\documentclass[draftclsnofoot]{IEEEtran}
\onecolumn
\usepackage{amsmath,amsfonts,amssymb}
\usepackage{dsfont}
\usepackage[T1]{fontenc}
\usepackage{graphicx}
\usepackage{subcaption}
\usepackage{cite}

\newtheorem{lemma}{Lemma}
\newtheorem{proposition}{Proposition}

\begin{document}
\ifCLASSINFOpdf
\else
\fi
\title{Gaussian integrals on symmetric spaces\\
(the complex case and beyond)}

\author{Salem Said (Laboratoire Jean Kuntzmann, Universit\'e Grenoble-Alpes)}

\maketitle

\section{Introduction}
The present work is concerned with Gaussian integrals,
$$
Z(\sigma) = \int_M\, e^{-d^{\hspace{0.02cm} 2}(x,x_o)/2\sigma^2}\hspace{0.02cm}\mathrm{vol}(dx) \hspace{1cm} \sigma > 0 \text{ and } x_o \in M
$$
where $M$ is a negatively-curved symmetric space, $d(x\hspace{0.02cm},x_o)$  the Riemannian distance and $\mathrm{vol}$~the Riemannian volume. There are two motivations for studying these integrals. The first one from statistics, through the study of Gaussian distributions on $M$~\cite{said2}. The second one comes from physics, in particular from the study of Chern-Simons theory~\cite{tierz}. 

The above Gaussian integrals may be viewed as perturbations of classical (Euclidean) Gaussian integrals,
$$
z(\sigma) = \int_{T_{x_o}M} e^{-\Vert v\Vert^2/2\sigma^2}\hspace{0.02cm}dv
$$
where the norm $\Vert v \Vert$ and volume measure $dv$ are those induced on $T_{x_o}M$ by the Riemannian metric of $M$.  Specifically,  one may show that (see Appendix \ref{app:perturbation})
$$
Z(\sigma) = z(\sigma)\left(1+ f_2\hspace{0.03cm}\sigma^2 + f_4\hspace{0.03cm}\sigma^4 + \ldots\hspace{0.03cm} \right)
$$
where, letting $d$ denote the dimension of $M$, the general coefficient $f_{2n}$ reads
$$
f_{2n} = \frac{1}{(2\pi)^{d/2}}\int_{T_{x_o}M}P_{\hspace{0.02cm}2n}(v)\hspace{0.02cm}e^{-\Vert v\Vert^2/2}\hspace{0.02cm}dv
$$
with $P_{\hspace{0.02cm}2n}$ an isotropy-invariant polynomial of degree $2n$ (here, isotropies are the isometries of $M$ that preserve $x_o$). The coefficients of $P_{\hspace{0.02cm}2n}$ can be written down in terms of the Riemann curvature tensor of $M$ at $x_o$ and $f_{2n}$ can then be evaluated using Wick's theorem, since it is just an expectation with respect to a normal probability density. 

Here, this perturbative approach will not be pursued. Rather, the aim is to directly evaluate the high-rank limit
$$
\lim_{r\rightarrow \infty}\hspace{0.02cm}\frac{1}{r^2}\log\!\left(\frac{Z(\sigma)}{z(\sigma)}\right) 
$$
where $r$ is the rank of $M$ and the limit is taken while $t = r\hspace{0.02cm}\sigma^2$ remains constant. The main idea is to adopt a two-step approach\,: in the complex case (the isometry group of $M$ is a complex Lie group) the high-rank limit admits a straightforward closed-form expression. In general, beyond the complex case, 
this limit can be obtained by solving a certain variational problem (minimising an energy functional over the space of probability distributions on $\mathbb{R}$).\hfill\linebreak Finally, this variational formulation can be exploited to recover the limit in closed-form, from the expression originally found in the complex case. 

Below, Section \ref{sec:background} fixes notations and provides some general background on the integrals $Z(\sigma)$ and $z(\sigma)$. Sections~\ref{sec:complex} and \ref{sec:hr1} deal with the complex case. In particular, Proposition \ref{prop:typeiv} gives a closed-form expression of~$Z(\sigma)$ (this generalises the expression of the Chern-Simons partition function from~\cite{marino}\cite{marinopaper}), while Proposition \ref{prop:free} uses this expression to evaluate the high-rank limit for each one of the four classical families of root systems. Section~\ref{sec:betanot2} moves beyond the complex case. Instead of four families of classical symmetric spaces, there are then eleven.\hfill\linebreak To mitigate this difficulty, 
Proposition \ref{prop:equilibrium} introduces the general variational formulation of the high-rank limit. Finally,~Section \ref{sec:cones} carries out the above-described two-step approach for two classes of symmetric spaces. Specifically, these are symmetric cones and classical symmetric domains.




















   
\section{Gaussian integrals} \label{sec:background}
Let $M$ be a Riemannian symmetric space, which is simply connected and has non-positive sectional curvatures. 
Assume, in addition, $M$ is not a Euclidean space. Consider the following integrals, here called Gaussian integrals,
\begin{equation} \label{eq:Z}
   Z(\sigma) = \int_M\, e^{-d^{\hspace{0.02cm} 2}(x,x_o)/2\sigma^2}\hspace{0.02cm}\mathrm{vol}(dx) \hspace{1cm} \sigma > 0 \text{ and } x_o \in M
\end{equation}
where $d(x\hspace{0.02cm},x_o)$ denotes Riemannian distance and $\mathrm{vol}$ Riemannian volume. Here, the point $x_o$ is arbitrary but fixed. The notation $Z(\sigma)$ is correct as the right-hand side of (\ref{eq:Z}) does not depend on $x_o\hspace{0.03cm}$. Indeed, the integral on the right-hand side is invariant under any isometry transformation of $M$. 

The behavior of $Z(\sigma)$ when $\sigma$ is close to zero is that of a classical (Euclidean) Gaussian integral. 
\begin{proposition} \label{prop:Ztaylor}
Let $Z(\sigma)$ be defined as in (\ref{eq:Z}). Then, 
\begin{equation} \label{eq:Ztaylor}
  Z(\sigma) = z(\sigma)\left(1+O\left(\sigma^2\right)\right)
\end{equation}
where $z(\sigma)$ is the classical Gaussian integral,
\begin{equation} \label{eq:z}
  z(\sigma) = \int_{T_{x_o}M} e^{-\Vert v\Vert^2/2\sigma^2}\hspace{0.02cm}dv
\end{equation}
Here, the norm $\Vert v \Vert$ and volume measure $dv$ are those induced on $T_{x_o}M$ by the Riemannian metric of $M$. Of~course, $z(\sigma) = (2\pi\sigma^2)^{d/2}$ where $d$ is the dimension of $M$.  
\end{proposition}
In certain cases, the classical integral $z(\sigma)$ is moreover related to the beautiful Macdonald-Mehta integral~\cite{etinghof} (Section 4). To make this connection, recall first the~integral formulas for the Cartan decomposition~\cite{helgasonbis} (Page 186),
\begin{align}
\label{eq:Zcartan} 
Z(\sigma) = \frac{\omega_M}{|W|}\int_{\mathfrak{a}} e^{-\Vert a\Vert^2/2\sigma^2}\hspace{0.02cm}\prod_{\lambda \in \Delta_+} |\sinh\lambda(a)|^{m_\lambda}\hspace{0.02cm}da \\[0.15cm]
\label{eq:zcartan}
z(\sigma) = \frac{\omega_M}{|W|}\int_{\mathfrak{a}} e^{-\Vert a\Vert^2/2\sigma^2}\hspace{0.02cm}\prod_{\lambda \in \Delta_+} |\lambda(a)|^{m_\lambda}\hspace{0.02cm}da
\end{align}
where the notation is the following~\cite{helgasonbis}\cite{helgason}. The symmetric space $M$ corresponds to a  symmetric pair~$(G,K)$. The~Lie~algebra $\mathfrak{g}$ of $G$ has a Cartan decomposition $\mathfrak{g} = \mathfrak{k} + \mathfrak{p}$ where $\mathfrak{k}$ is the Lie algebra of $K$ and $\mathfrak{p} \simeq T_{x_o}M$. If $\mathfrak{a}$ is a maximal abelian subspace of $\mathfrak{p}$, then the corresponding set of positive restricted roots is denoted by $\Delta_+\hspace{0.03cm}$, and the dimension of the root subspace corresponding to $\lambda \in \Delta_+$ is denoted by $m_\lambda\hspace{0.03cm}$. Moreover, $|W|$ is the order of the Weyl group $W$, the group of orthogonal transformations of $\mathfrak{a}$, generated by reflections in the root hyperplanes. Finally, $\omega_M$ is a certain constant, which only depends on the space $M$.
   
To bring in the Macdonald-Mehta integral, assume that $m_\lambda = 2$ does not depend on $\lambda$. Then~\cite{etinghof} (Section 4)
\begin{equation} \label{eq:macdonald1}
\int_{\mathfrak{a}} e^{-\Vert a\Vert^2/2}\hspace{0.02cm}\prod_{\lambda \in \Delta_+} |\lambda(a)|^2\hspace{0.02cm}da = (2\pi)^{r/2}\prod_{\lambda \in \Delta_+} \left(\Vert \lambda\Vert^2\middle/2\right)\prod^r_{j=1}\Gamma(1+d_j)
\end{equation}
where $r$ denotes the dimension of $\mathfrak{a}$, $\Vert \lambda \Vert$ is the norm of $\lambda \in \mathfrak{a}^*$ with respect to the Riemannian metric of $M$, and~$d_j$ are the degrees of the generators of the algera of $W$-invariant polynomial functions on $\mathfrak{a}$~\cite{humphreys}. By a straightforward change of variables, (\ref{eq:zcartan}) now yields
\begin{equation} \label{eq:macdonald2}
z(\sigma) = \frac{\omega_M}{|W|}\hspace{0.05cm}\sigma^d \times (2\pi)^{r/2} \prod_{\lambda \in \Delta_+} \left(\Vert \lambda\Vert^2\middle/2\right)\prod^r_{j=1}\Gamma(1+d_j)
\end{equation}
The general expression of the Macdonald-Mehta integral is in no way limited to the special case $m_\lambda = 2$~\cite{etinghof}\cite{opdam}. Still, this is enough when the isometry group $G$ is assumed to have a complex structure, as in the following section.     

\section{The complex case} \label{sec:complex}
The Gaussian integral (\ref{eq:Z}) admits a closed-form expression whenever $M$ is a so-called type IV symmetric space. This means the Lie algebra $\mathfrak{g}$ of $G$ is the realification of a simple complex Lie algebra~\cite{helgason} (Chapter X, Page 516). 
The following proposition not only expresses (\ref{eq:Z}), but even the more general integrals
\begin{equation} \label{eq:Zl}
  Z(\sigma,\tau) = \int_M\,e^{-d^{\hspace{0.02cm} 2}(x,x_o)/2\sigma^2}\hspace{0.02cm}\Phi_\tau(x)\hspace{0.02cm}\mathrm{vol}(dx)
\end{equation}
where $\Phi_\tau$ denotes the spherical function corresponding to $\tau \in \mathfrak{a}^*_{\scriptscriptstyle\, \mathbb{C}}$~\cite{helgasonbis} (Page 418). This is an eigenfunction of the Riemannian Laplace operator of $M$, with eigenvalue $(\tau,\tau) - (\rho,\rho)$ where $\rho$ is the sum of $\lambda \in \Delta_+$ (the brackets denote the bilinear form induced by the Riemannian metric). In particular, the original $Z(\sigma)$ from (\ref{eq:Z}) is  just $Z(\sigma,\rho)$.
\begin{proposition}  \label{prop:typeiv}
    If $M$ is a type IV symmetric space, then the integrals (\ref{eq:Zl}) are given by
\begin{equation} \label{eq:typeiv}
  Z(\sigma,\tau) = \omega_M\hspace{0.03cm}\sigma^d\times (2\pi)^{r/2}\hspace{0.03cm}\varpi(\rho/2)\times e^{\frac{\sigma^2}{2}\left((\tau,\tau)+(\rho,\rho)\right)}\prod_{\lambda \in \Delta_+}\mathrm{sch}\left(\sigma^2(\lambda,\tau)\right)
\end{equation}
In particular, for the integral (\ref{eq:Z}),
\begin{equation} \label{eq:typeivZ}
  Z(\sigma) = \frac{\omega_M}{2^{|\Delta_+|}}\times (2\pi\sigma^2)^{r/2} \times e^{\sigma^2\hspace{0.02cm}(\rho,\rho)}\prod_{\lambda \in \Delta_+} \sinh\left(\sigma^2(\lambda,\rho)\right)
\end{equation}
Here, $\varpi(\tau) = \prod_{\lambda \in \Delta_+} (\lambda,\tau)$ for any $\tau \in \mathfrak{a}^*_{\scriptscriptstyle\, \mathbb{C}}\hspace{0.03cm}$ and $\mathrm{sch}(x) = \sinh(x)/x$. 
\end{proposition}
Propositions \ref{prop:Ztaylor} and \ref{prop:typeiv} can be applied jointly, leading to some amusing results. For example (see Paragraph \ref{subsec:prooftaylor}), they~can be used to recover the well-known strange formula of Freudenthal and de Vries~\cite{freudenthal}. For another application, note~that replacing $\tau = \rho$ into (\ref{eq:typeiv}) yields
\begin{equation} \label{eq:typeivZ_bis}
Z(\sigma)  = 
\omega_M\hspace{0.03cm}\sigma^d\times (2\pi)^{r/2}\hspace{0.03cm}\varpi(\rho/2)\times e^{\sigma^2(\rho,\rho)}\prod_{\lambda \in \Delta_+}\mathrm{sch}\left(\sigma^2(\lambda,\rho)\right) 
\end{equation}
However, it then follows from (\ref{eq:Ztaylor}) that
\begin{equation} \label{eq:typeivz_bis}
  z(\sigma) = 
\omega_M\hspace{0.03cm}\sigma^d\times (2\pi)^{r/2}\hspace{0.03cm}\varpi(\rho/2)
\end{equation}
This can be used to prove (\ref{eq:macdonald2}), and therefore to recover (\ref{eq:macdonald1}). To do so, it enough to replace into (\ref{eq:typeivz_bis}) the identity~\cite{macdonald}
\begin{equation} \label{eq:magicformula}
  \check{\varpi}(\rho/2) = 
\prod^r_{j=1}\Gamma(d_j)
\end{equation}   
where $\check{\varpi}(\tau) = \prod_{\lambda \in \Delta_+} (\check{\lambda},\tau)$ in terms of the coroots $\check{\lambda} = 2\lambda/\Vert\lambda\Vert^2$, and recall $|W| = \prod^r_{j=1}d_j$ (see~\cite{humphreys}, Page 62).     

\section{The high-rank limit} \label{sec:hr1}
The classical type IV symmetric spaces correspond to the four classical families of reduced root systems~\cite{helgason}. 
Let~$M$ be a classical type IV symmetric space and denote $r$ its rank (so that $r$ is the dimension of  $\mathfrak{a}$, as above). Then, consider the limit
\begin{equation} \label{eq:tilde_free}
  F(t) = \lim_{r \rightarrow \infty}\frac{1}{r^2}\log\!\left(Z(\sigma)\middle/z(\sigma)\right) 
\end{equation}
which is taken as the rank $r$ increases to infinity while $t = r\hspace{0.02cm}\sigma^2$ remains constant. 

Recall that the Riemannian metric of $M$ is unique up to a multiplicative constant. In particular, for any $\lambda \in \Delta_+\hspace{0.03cm}$, $\Vert \lambda\Vert^2 = c, 2c$ or $4c$ where $c > 0$. For the limit (\ref{eq:tilde_free}) to exist, it is enough that $c$ should not depend on the rank $r$. Here, without loss of generality, $c$ is chosen equal to $1/4$ (the sectional curvatures of $M$ are then always $\geq -1$).

\begin{proposition} \label{prop:free}
Let~$M$ be a classical type IV symmetric space, and assume its Riemannian metric is chosen so 
$\Vert \lambda\Vert^2 = 1/4, 1/2$ or $1$ for any $\lambda \in \Delta_+\hspace{0.02cm}$. Then, the limit in (\ref{eq:tilde_free}) exists for any $t > 0$ and has the following expression
\begin{itemize}
 \item[--] $A_{r}$ root system\,: $M = \mathrm{SL}(r+1,\mathbb{C})/\mathrm{SU}(r+1)$, and 
\begin{equation} \label{eq:asymptotics_A}
  F(t) = \frac{t}{3} + \int^1_0 (1-x)\hspace{0.02cm}\log\left(\mathrm{sch}\left(tx/2\right)\right)dx
\end{equation}
\item[--] $B_{r}$ root system\,: $M = \mathrm{SO}(2r+1,\mathbb{C})/\mathrm{SO}(2r+1)$, and
\begin{equation} \label{eq:asymptotics_B}
F(t) = \frac{t}{3}  + \frac{1}{2}\int^2_0 (2-x)\hspace{0.02cm}\log\left(\mathrm{sch}\left(tx/2\right)\right)dx
\end{equation}
\item[--] $C_{r}$ root system\,: $M = \mathrm{Sp}(r,\mathbb{C})/\mathrm{Sp}(r)$, and
\begin{equation} \label{eq:asymptotics_C}
F(t) = \frac{t}{3}  + \frac{1}{2}\int^2_0 (2-x)\hspace{0.02cm}\log\left(\mathrm{sch}\left(tx/2\right)\right)dx
\end{equation}
\item[--] $D_{r}$ root system\,: $M = \mathrm{SO}(2r,\mathbb{C})/\mathrm{SO}(2r)$, and
\begin{equation} \label{eq:asymptotics_D}
F(t) = \frac{t}{3}  + \frac{1}{2}\int^2_0 (2-x)\hspace{0.02cm}\log\left(\mathrm{sch}\left(tx/2\right)\right)dx
\end{equation}
\end{itemize}
\end{proposition}
The above proposition shows that the limit $F(t)$ is the same for all the root systems other than $A_r$ (last three cases). 
From the proof of this proposition, it will also become clear that this limit holds with the following error estimate
\begin{equation} \label{eq:error1/r}
\frac{1}{r^2}\log\!\left(Z(\sigma)\middle/z(\sigma)\right) = F(t) + O\!\left(\frac{1}{r}\right)
\end{equation}
where the $O(1/r)$ bound holds uniformly in $t$ whenever $t$ belongs to a bounded interval. 

The integrals in Proposition \ref{prop:free} can be evaluated explicitly, using the two identities
\begin{equation} \label{eq:dilog}
\int^a_{0} \log\left(1 - e^{-tx}\right)\hspace{0.02cm}dx = \frac{\mathrm{Li}_{\hspace{0.02cm}2}(e^{-ta}) - \mathrm{Li}_{\hspace{0.02cm}2}(1)}{t} \hspace{2.4cm}
\end{equation}
\begin{equation} \label{eq:trilog}
\int^a_{0} x\log\left(1 - e^{-tx}\right)\hspace{0.02cm}dx = \frac{a}{t}\mathrm{Li}_{\hspace{0.02cm}2}(e^{-ta}) + 
\frac{\mathrm{Li}_{\hspace{0.02cm}3}(e^{-ta}) - \mathrm{Li}_{\hspace{0.02cm}3}(1)}{t^2} 
\end{equation}
\noindent where $\mathrm{Li}_{\hspace{0.02cm}2}$ and $\mathrm{Li}_{\hspace{0.02cm}3}$ denote the dilogarithm and trilogarithm functions~\cite{irresistible}.

\section{Beyond the complex case} \label{sec:betanot2}   
Here, the assumption that $M$ should be a type IV symmetric space (as in Sections \ref{sec:complex} and \ref{sec:hr1}) will be dropped. Rather, $M$ is allowed to be any classical irreducible symmetric space of non-compact type. Instead of four families, as in Section \ref{sec:hr1}, $M$ now belongs to one of eleven families (these are the previous four and seven new ones~\cite{helgason}). 

The object of interest will be the integral which appears in (\ref{eq:Zcartan}),
\begin{equation} \label{eq:generalI}
  Z_\mathfrak{a}(\sigma) = \int_{\mathfrak{a}} e^{-\Vert a\Vert^2/2\sigma^2}\hspace{0.02cm}\prod_{\lambda \in \Delta_+} |\sinh\lambda(a)|^{m_\lambda}\hspace{0.02cm}da
\end{equation}
with the aim of understanding the behavior of $(1/r^2)\log\!\left(Z_\mathfrak{a}(\sigma)\right)$ when $r \rightarrow \infty$ while $t = r\hspace{0.02cm}\sigma^2$ remains constant ($r$~denotes the rank of $M$, just as in Section \ref{sec:hr1}). There are now five (instead of four) root systems, which yield the following cases in (\ref{eq:generalI}) -- see~\cite{helgason} (Chapter X, Table VI)~\cite{knapp} (Page 424)\,:
\begin{itemize}
 \item[--] $A_{r}$ root system\,:
\begin{equation} \label{eq:IA}
  Z_\mathfrak{a}(\sigma) = \int_{\mathfrak{a}} e^{-\frac{c}{\mathstrut 2\sigma^{\mathstrut\scriptscriptstyle 2}}\sum^{r+1}_{i=1}a^2_i}\hspace{0.02cm}\prod_{i<j} |\sinh(a_i-a_j)|^{\beta}\hspace{0.02cm}da
\end{equation}
 \item[--] $B_{r}$ root system\,:
\begin{equation} \label{eq:IB}
  Z_\mathfrak{a}(\sigma) = \int_{\mathfrak{a}} e^{-\frac{c}{\mathstrut 2\sigma^{\mathstrut\scriptscriptstyle 2}}\sum^r_{i=1}a^2_i}\hspace{0.02cm}\prod_{i<j} |\sinh(a_i-a_j)\sinh(a_i+a_j)|^{\beta}\prod_i|\sinh(a_i)|^{\alpha}\hspace{0.02cm}da
\end{equation}
 \item[--] $C_{r}$ root system\,:
\begin{equation} \label{eq:IC}
  Z_\mathfrak{a}(\sigma) = \int_{\mathfrak{a}} e^{-\frac{c}{\mathstrut 2\sigma^{\mathstrut\scriptscriptstyle 2}}\sum^r_{i=1}a^2_i}\hspace{0.02cm}\prod_{i<j} |\sinh(a_i-a_j)\sinh(a_i+a_j)|^{\beta}\prod_i|\sinh(2a_i)|^{\gamma}\hspace{0.02cm}da
\end{equation}
 \item[--] $D_{r}$ root system\,:
\begin{equation} \label{eq:ID}
  Z_\mathfrak{a}(\sigma) = \int_{\mathfrak{a}} e^{-\frac{c}{\mathstrut 2\sigma^{\mathstrut\scriptscriptstyle 2}}\sum^r_{i=1}a^2_i}\hspace{0.02cm}\prod_{i<j} |\sinh(a_i-a_j)\sinh(a_i+a_j)|^{\beta}\hspace{0.02cm}da
\end{equation}
 \item[--] $(BC)_{r}$ root system\,:
\begin{equation} \label{eq:IBC}
  Z_\mathfrak{a}(\sigma) = \int_{\mathfrak{a}} e^{-\frac{c}{\mathstrut 2\sigma^{\mathstrut\scriptscriptstyle 2}}\sum^r_{i=1}a^2_i}\hspace{0.02cm}\prod_{i<j} |\sinh(a_i-a_j)\sinh(a_i+a_j)|^{\beta}\prod_i|\sinh(a_i)|^{\alpha}
|\sinh(2a_i)|^{\gamma}\hspace{0.02cm}da
\end{equation}
\end{itemize}
Here, $c > 0$ is an arbitrary constant which determines the choice of the Riemannian metric of $M$, and each 
$a \in \mathfrak{a}$ is identified with an array of real components $(a_i)$ (with the additional requirement $\sum a_i = 0$ for $A_r$ root systems). The values of the multiplicities $\alpha,\beta$ and $\gamma$ are recalled in Table \ref{tab:classif} below. The notation used in this table (specifically, the naming of classical Lie groups in the first two columns) is the standard notation found in~\cite{helgason}\cite{knapp}. 

The desired asymptotic behavior of $\log\!\left(Z_\mathfrak{a}(\sigma)\right)$ is given by Proposition \ref{prop:equilibrium}, which involves the following notation. Let $\mu$ be a probability measure on $\mathbb{R}$ and consider the energy functionals (compare to~\cite{deift}, Chapter 6)
\begin{align} 
\label{eq:functionala}
I_{A}(\mu) = &\frac{c}{2t} \int\, x^2\hspace{0.02cm}\mu(dx) - \frac{\beta}{2} \int\int\,\ell(x-y)\hspace{0.02cm}\mu(dx)\mu(dy) \\
\label{eq:functionalc} I_{C}(\mu) = & \frac{c}{2t} \int\, x^2\hspace{0.02cm}\mu(dx) - \frac{\beta}{2} \int\int\left(\ell(x-y)+\ell(x+y)\right)\mu(dx)\mu(dy) \\
\label{eq:functionalbc} I_{BC}(\mu) = & \int\,\left((c/2t)\hspace{0.02cm}x^2 - \beta(\delta-1)\hspace{0.02cm}\ell(x)\right)\mu(dx) - \frac{\beta}{2} \int\int\left(\ell(x-y)+\ell(x+y)\right)\mu(dx)\mu(dy) 
\end{align} 
where $\ell(x) = \log(|\sinh(x)|)$ and $\delta \geq 1$ a constant to be defined in the proposition. In addition,  let $E_{A} = \inf I_{A}(\mu)$, $E_{C} = \inf I_{C}(\mu)$ and $E_{BC} = \inf I_{BC}(\mu)$ (where the $\inf$ is over probability measures $\mu$ on $\mathbb{R}$).
\begin{proposition} \label{prop:equilibrium}
The ``equilibrium energies" $E_A$, $E_C\hspace{0.01cm}$, and $E_{BC}$ (the infima just defined) are finite.
In addition, if the constant $c$ (which was introduced in (\ref{eq:IA})-(\ref{eq:IBC})) does not depend on $r$, then the following limits hold when $r \rightarrow \infty$ while $t = r\hspace{0.02cm}\sigma^2$ remains constant (this is true for any $t > 0$)
\begin{itemize}
\item[--] $A_{r}$ root system\,:  
\begin{equation} \label{eq:lima}
\lim_{r\rightarrow \infty}\frac{1}{r^2}\log\!\left(Z_\mathfrak{a}(\sigma)\right) = -E_{A} 
\end{equation}
where $Z_\mathfrak{a}(\sigma)$ is given by (\ref{eq:IA}).
\item[--] $B_{r}$ root system\,:    
\begin{equation} \label{eq:limb}
\lim_{r\rightarrow \infty}\frac{1}{r^2}\log\!\left(Z_\mathfrak{a}(\sigma)\right) = -E_{BC}
\end{equation}
where $Z_\mathfrak{a}(\sigma)$ is given by (\ref{eq:IB}), under the additional assumption that $\alpha = \beta(\delta - 1)\hspace{0.02cm}r$ for some constant $\delta \geq 1$. 
\item[--] $C_r$ or $D_r$ root system\,:
\begin{equation} \label{eq:limcd}
\lim_{r\rightarrow \infty}\frac{1}{r^2}\log\!\left(Z_\mathfrak{a}(\sigma)\right) = -E_{C}
\end{equation}
where $Z_\mathfrak{a}(\sigma)$ is given by (\ref{eq:IC}) or (\ref{eq:ID}). 
\item[--] $(BC)_r$ root system\,:
\begin{equation} \label{eq:limbc}
\lim_{r\rightarrow \infty}\frac{1}{r^2}\log\!\left(Z_\mathfrak{a}(\sigma)\right) = -E_{BC}
\end{equation}
where $Z_\mathfrak{a}(\sigma)$ is given by (\ref{eq:IBC}), under the additional assumption that $\alpha = \beta(\delta - 1)\hspace{0.02cm}r$ for some constant $\delta \geq 1$. 
On~the other hand, the limit is equal to $-E_{C}$ if $\alpha = 4$ or $0$ independent of $r$ (as in the last line of Table \ref{tab:classif}). 
\end{itemize}
\end{proposition}
Note that the limits for $B_{r}$ and $(BC)_{r}$ root systems require the additional assumption that $\alpha = \beta (\delta - 1)\hspace{0.02cm}r$. In~the~notation of Table \ref{tab:classif}, this happens when $r = p$ while $\alpha = \beta(q-p)$ and $q$ is assumed to be equal to $\delta\hspace{0.02cm}p$. If~$\delta = 1$ (that is, $q = p$) then $E_{BC} = E_C\hspace{0.02cm}$, as one may see from (\ref{eq:functionalbc}). 

\begin{table}[t!]
\begin{center}
\begin{tabular}{|c|c|c|c|c|c|c|c|c|}
\hline 
$G$ & $K$ & root system & rank $r$ & dimension $d$ & $\alpha$ & $\beta$ &  $\gamma$ & remark \\[0.1cm]
\hline
$\mathrm{SL}(r,\mathbb{R})$ & $\mathrm{SO}(r)$ & $A$ & $r-1$ &  $(r-1)(r+2)/2$ &      & $1$ & &\\[0.1cm] 
$\mathrm{SL}(r,\mathbb{C})$ & $\mathrm{SU}(r)$ & $A$ & $r-1$ &  $(r-1)(r+1)$        &      & $2$ & & complex case \\[0.1cm] 
$\mathrm{SL}(r,\mathbb{H})$ & $\mathrm{Sp}(r)$ & $A$ & $r-1$ &  $ (r-1)(2r+1)$        &      & $4$ & & \\[0.1cm] 
$\mathrm{SO}_{o}(p,q)$               & $SO(p) \times SO(q)$ & $B$ or $D$ & $p$   & $p\hspace{0.02cm}q$ &    $q-p$   & $1$  &  & $p \leq q$ \\[0.1cm]
$\mathrm{SO}(2r+1,\mathbb{C})$               & $SO(2r+1)$ & $B$ & $r$   & $r(2r+1)$ &    $2$   & $2$  &  &  complex case \\[0.1cm]
$\mathrm{SO}(2r,\mathbb{C})$               & $SO(2r)$ & $D$ & $r$   & $r(2r-1)$ &       & $2$  &  &  complex case \\[0.1cm]
$\mathrm{Sp}(r,\mathbb{R})$ & $\mathrm{U}(r)$ & $C$ & $r$ & $r(r+1)$ & & $1$ & $1$ & \\[0.1cm]
$\mathrm{Sp}(r,\mathbb{C})$ & $\mathrm{Sp}(r)$ & $C$ & $r$ &  $r(2r+1)$ &  & $2$ & $2$ & complex case \\[0.1cm]
$\mathrm{SU}(p,q)$ & $\mathrm{S}(\mathrm{U}(p) \times \mathrm{U}(q))$  & $BC$ or $C$ & $p$ & $2p\hspace{0.02cm}q$ & $2(q-p)$ &
$2$ & $1$ & $p \leq q$ \\[0.1cm]
$\mathrm{Sp}(p,q)$ & $\mathrm{Sp}(p) \times \mathrm{Sp}(q)$ & $BC$ or $C$ & $p$ & $4p\hspace{0.02cm}q$ & $4(q-p)$ & $4$ & $3$ & $p \leq q$ \\[0.1cm]
$\mathrm{SO}^*(2n)$ & $\mathrm{U}(n)$ & $BC$ or $C$ & $\lfloor n/2 \rfloor$ & $n(n-1)$ & $4$ or $0$ & 4 & 1 & $\alpha = 4$ if $n$ is odd \\[0.2cm]
\hline
\end{tabular}
\end{center}
\caption{Classification of the classical irreducible symmetric spaces of non-compact type}
\label{tab:classif}
\end{table}







\section{Application to cones and domains} \label{sec:cones}
The aim here is to jointly apply Propositions \ref{prop:free} and \ref{prop:equilibrium} to certain special classes of symmetric spaces.

\subsection{Symmetric cones}
For $\beta = 1, 2, 4$, consider the symmetric space $M$ determined by the symmetric pair $(G,K)$, 
$$
\begin{array}{ccc}
\beta & G & K \\[0.1cm]
1      & \mathrm{GL}(r,\mathbb{R}) & \mathrm{O}(r) \\[0.1cm]
2      & \mathrm{GL}(r,\mathbb{C}) & \mathrm{U}(r) \\[0.1cm]
4      & \mathrm{GL}(r,\mathbb{H}) & \mathrm{Sp}(r)
\end{array}
$$
Then, $M$ is a symmetric cone in the sens of~\cite{faraut}. In fact, it is the cone of positive-definite matrices over $\mathbb{R}$, $\mathbb{C}$, $\mathbb{H}$, for $\beta = 1, 2, 4$ (respectively). For this $M$, the integral (\ref{eq:Zcartan}) takes on the  form
\begin{equation} \label{eq:Zcone}
  Z(\sigma) = \frac{\omega_\beta}{r!}\int_{\mathbb{R}^r} e^{-\Vert a\Vert^2/2\sigma^2}\hspace{0.02cm}\prod_{i<j} |\sinh(a_i-a_j)|^\beta\hspace{0.02cm}da \\[0.15cm]
\end{equation}
On the other hand, 
the~integral (\ref{eq:zcartan}) becomes
\begin{equation} \label{eq:zcone}
  z(\sigma) = \frac{\omega_\beta}{r!}\int_{\mathbb{R}^r} e^{-\Vert a\Vert^2/2\sigma^2}\hspace{0.02cm}\prod_{i<j} |a_i-a_j|^\beta\hspace{0.02cm}da \\[0.15cm]
\end{equation}
where the constant $\omega_\beta$ also depends on $r$. The Riemannian metric of $M$ will be determined by $\Vert a\Vert^2 = 4 \sum a^2_j\hspace{0.03cm}$. 

Let $Z_\mathfrak{a}(\sigma)$ and $z_\mathfrak{a}(\sigma)$ denote $Z(\sigma)$ and $z(\sigma)$ with the factor $\omega_\beta/r!$ dropped from (\ref{eq:Zcone}) and (\ref{eq:zcone}). Proposition \ref{prop:equilibrium} 
yields the following asymptotic behavior of $\log\!\left(Z_\mathfrak{a}(\sigma)\right)$\,: let $E_A$ be defined as in Proposition \ref{prop:equilibrium} and write $E_A = E_\beta(t)$, to~highlight dependence on $\beta$ and $t$. With this notation, it is possible to show that Proposition \ref{prop:equilibrium} implies the limit
\begin{equation} \label{eq:Zcone_asymp}
  \lim_{r\rightarrow \infty}\frac{1}{r^2}\log\!\left(Z_\mathfrak{a}(\sigma)\right) = -\frac{\beta}{2}\hspace{0.03cm}
E_{2}\!\left(\frac{\beta}{2}t\right)
\end{equation}
where the limit is taken while $t = r\hspace{0.02cm}\sigma^2$ remains fixed. For the asymptotic behavior of $\log\!\left(z_\mathfrak{a}(\sigma)\right)$, note first that, using the Macdonald-Mehta integral (in its general form~\cite{etinghof}\cite{opdam}),
\begin{equation} \label{eq:zcone_mehta}
   z_\mathfrak{a}(\sigma) = \left(\frac{\sigma}{2}\right)^{\!d_\beta}\!\times (8\pi)^{r/2}\prod^r_{j=1} \frac{\Gamma(1+j\beta/2)}{\Gamma(1+\beta/2)}
\end{equation}
where $d_\beta$ is the dimension of $M$ (this is $d_\beta = r + r(r-1)\beta/2$). Starting from this identity, it is possible to show
\begin{equation} \label{eq:zcone_asymp}
  \lim_{r\rightarrow \infty}\frac{1}{r^2}\log\!\left(z_\mathfrak{a}(\sigma)\right) = -\frac{\beta}{2}\hspace{0.03cm}
e_{2}\!\left(\frac{\beta}{2}t\right) \hspace{1cm} e_2(t) = \frac{3}{4} - \frac{1}{2}\log(t/4)
\end{equation}
where the limit is taken while $t = r\hspace{0.02cm}\sigma^2$ remains fixed. The proofs of (\ref{eq:Zcone_asymp}) and (\ref{eq:zcone_asymp}) are provided in Paragraph \ref{app:cones}. 
Subtracting (\ref{eq:zcone_asymp}) from (\ref{eq:Zcone_asymp}), it is clear that
$$
  \lim_{r \rightarrow \infty}\frac{1}{r^2}\log\!\left(Z(\sigma)\middle/z(\sigma)\right) = -\frac{\beta}{2}\left[
E_{2}\!\left(\frac{\beta}{2}t\right) -e_{2}\!\left(\frac{\beta}{2}t\right)\right]
$$
However (see details in Paragraph \ref{app:cones}), when $\beta = 2$, this limit is given by (\ref{eq:asymptotics_A}). Therefore, if $F(t)$, $Z(\sigma)$ and~$z(\sigma)$ are given by (\ref{eq:asymptotics_A}), (\ref{eq:Zcone}) and (\ref{eq:zcone}),
\begin{equation}\label{eq:cone_ratio_asymp}
  \lim_{r \rightarrow \infty}\frac{1}{r^2}\log\!\left(Z(\sigma)\middle/z(\sigma)\right) = \frac{\beta}{2}\hspace{0.03cm}F\!\left(\frac{\beta}{2}t\right)
\end{equation}
where the limit is taken while $t = r\hspace{0.02cm}\sigma^2$ remains fixed.


\subsection{Classical domains}
For $\beta = 1, 2, 4$, consider the symmetric space $M$ determined by the symmetric pair $(G,K)$, 
$$
\begin{array}{ccc}
\beta & G & K \\[0.1cm]
1      & \mathrm{Sp}(r,\mathbb{R}) & \mathrm{U}(r) \\[0.1cm]
2      & \mathrm{SU}(r,r) & \mathrm{S}(\mathrm{U}(r) \times \mathrm{U}(r)) \\[0.1cm]
4      & \mathrm{SO}^*(4r) & \mathrm{U}(2r)
\end{array}
$$
Then, $M$ can be realised as one of the classical symmetric domains equipped with its Bergman metric~\cite{piatetski}. 
For~example, when~$\beta = 2$, this is the unit disc, with respect to the Euclidean operator norm, in the space of $r \times r$ complex matrices.

Looking up Table \ref{tab:classif}, the integral (\ref{eq:generalI}) is here seen to be of the form (\ref{eq:IC}) with $\gamma = 1$. To apply Proposition \ref{prop:equilibrium}, write $E_C = E_\beta(t)$, and note that (\ref{eq:functionalc}) implies $E_\beta(t) = (\beta/2)E_2(\beta t/2)$.  Therefore (just as in (\ref{eq:Zcone_asymp})),
\begin{equation} \label{eq:Zdomain_asymp}
  \lim_{r\rightarrow \infty}\frac{1}{r^2}\log\!\left(Z_\mathfrak{a}(\sigma)\right) = -\frac{\beta}{2}\hspace{0.03cm}
E_{2}\!\left(\frac{\beta}{2}t\right)
\end{equation}
where the limit is taken while $t = r\hspace{0.02cm}\sigma^2$ remains fixed.

To identify the function $E_2(t)$, it is helpful to resort back to the complex case, and consider $M = \mathrm{Sp}(r,\mathbb{C})/\mathrm{Sp}(r)$ (see (\ref{eq:asymptotics_C}), Proposition \ref{prop:free}). For this $M$, the integral (\ref{eq:generalI}) is also of the form (\ref{eq:IC}), with $\beta = 2$ and $\gamma = 2$ (Table \ref{tab:classif}).\hfill\linebreak 
Application of Proposition \ref{prop:equilibrium} then yields the same limit $-E_C = -E_2(t)$ (by~(\ref{eq:functionalc}), this~does not depend on $\gamma$). Keeping this in mind, $E_2(t)$ can be recovered
from Proposition \ref{prop:free}\,: if $F(t)$ is given by (\ref{eq:asymptotics_C}), Proposition \ref{prop:free} states
\begin{equation} \label{eq:free_symplectic}
F(t) = \lim_{r \rightarrow \infty}\frac{1}{r^2}\log\!\left(Z(\sigma)\middle/z(\sigma)\right)  
\hspace{1cm} \text{($t = r\sigma^2$ fixed)}
\end{equation}
where the integrals $Z(\sigma)$ and $z(\sigma)$ are those of $M = \mathrm{Sp}(r,\mathbb{C})/\mathrm{Sp}(r)$. From (\ref{eq:Zcartan}) and (\ref{eq:zcartan}), it is quite clear that $\left.Z(\sigma)\middle/z(\sigma)\right. = \left.Z_\mathfrak{a}(\sigma)\middle/z_\mathfrak{a}(\sigma)\right.$. However, as just explained, $(1/r^2)\log\!\left(Z_\mathfrak{a}(\sigma)\right)$ converges to the limit $-E_2(t)$. On~the~other hand, 
\begin{equation} \label{eq:za_symplectic}
  z_\mathfrak{a}(\sigma) = 
  \int_{\mathbb{R}^r} e^{-\Vert a\Vert^2/2\sigma^2}\hspace{0.02cm}\prod_{i<j} |(a_i-a_j)(a_i+a_j)|^2\prod_i|2a_i|^2\hspace{0.02cm}da
\end{equation}
which is given by the Macdonald-Mehta formula (\ref{eq:macdonald2}) (with $d_j = 2j$ as found in~\cite{humphreys}, Page 59),
\begin{equation} \label{eq:symplectic_macdonald}
  z_\mathfrak{a}(\sigma) = 
  (\sigma^2/4)^{r^2}\times (8\pi\sigma^2)^{r/2}\prod^r_{j=1}\Gamma(1+2j)
\end{equation}
From this expression, it can be shown (the proof is ommitted as it closely mirrors the proof of (\ref{eq:zcone_asymp}))
\begin{equation} \label{eq:z_symplectic_asymp}
    \lim_{r\rightarrow \infty}\frac{1}{r^2}\log\!\left(z_\mathfrak{a}(\sigma)\right) = 
-e_{2}\!\left(t\right) \hspace{1cm} e_2(t) = \frac{3}{2} - \log(t/2)
\end{equation}
where the limit is taken while $t = r\hspace{0.02cm}\sigma^2$ remains fixed. Putting together all of the above,
$$
  F(t) = -(E_2(t)-e_2(t)) 
$$
which means that
\begin{equation} \label{eq:asymptotic_domains_conclusion}
  E_2(t) = -(F(t)-e_2(t))
\end{equation}
where $F(t)$ is given by (\ref{eq:asymptotics_C}) and $e_2(t)$ is given by (\ref{eq:z_symplectic_asymp}). 

This is an analytic expression of $E_2(t)$ which also yields all of the limits in (\ref{eq:Zdomain_asymp}).  

\section{Proofs of the main results}

\subsection{Proof of Proposition \ref{prop:Ztaylor}} \label{subsec:prooftaylor}
To begin from the end, (\ref{eq:z}) is a classical Gaussian integral over the $d$-dimensional Euclidean space $T_{x_o}M$. Therefore, $z(\sigma) = (2\pi\sigma^2)^{d/2}$. To prove (\ref{eq:Ztaylor}), note that the Riemannian exponential map is a diffeomorphism from $T_{x_o}M$ to $M$. Denoting this by $\mathrm{Exp}$, one has
\begin{equation} \label{eq:prooftaylor1}
Z(\sigma) = \int_{T_{x_o}M}e^{-\Vert v\Vert^2/2\sigma^2}\hspace{0.02cm}\mathrm{Exp}^*(\mathrm{vol})(dv)
\end{equation}
by pulling the integral (\ref{eq:Z}) back to $T_{x_o}M$. However~\cite{helgason} (Page 215), using the notation introduced after (\ref{eq:Zcartan}) and (\ref{eq:zcartan}),
\begin{equation} \label{eq:prooftaylor2}
\mathrm{Exp}^*(\mathrm{vol})(dv) = \det\left(\sum^\infty_{n=0} \frac{\mathrm{ad}^{2n}(v)}{(2n+1)!} \right)dv
\end{equation}
where $\mathrm{ad}^2(v) : \mathfrak{p} \rightarrow \mathfrak{p}$ with $\mathrm{ad}^2(v)\hspace{0.02cm}u =[v,[v,u]]$ for $v,u \in \mathfrak{p} \simeq T_{x_o}M$ (the bracket is that of the Lie algebra $\mathfrak{g})$. 
Now, $\mathrm{ad}^2(v)$ is always diagonalisable. If $v = \mathrm{Ad}(k)\hspace{0.02cm}a$ with $k \in K$ and $a\in \mathfrak{a}$, the eigenvalues are either $0$ or $\lambda^2(a)$ with multiplicity $m_\lambda$ (for each $\lambda \in \Delta_+$).  The above determinant is therefore equal to 
\begin{align*}
\prod_{\lambda \in \Delta_+} \left(\mathrm{sch}\left(\lambda(a)\right)\right)^{m_\lambda} & = 1 + \frac{1}{6}\sum_{\lambda \in \Delta_+} m_\lambda\hspace{0.03cm}\lambda^2(a) + O\left(\Vert a \Vert^4\right) \\[0.1cm]
&= 1 + \frac{1}{6}\,\mathrm{tr}\left(\mathrm{ad}^2(v)\right) + O\left(\Vert v \Vert^4\right)
\end{align*}
where $\mathrm{sch}(x) = \sinh(x)/x$ and $\mathrm{tr}$ denotes the trace. Replacing into (\ref{eq:prooftaylor1}) and (\ref{eq:prooftaylor2}), 
$$
Z(\sigma) = z(\sigma)\,\int_{T_{x_o}M}G(v\hspace{0.03cm};\sigma)\left(1 + \frac{1}{6}\mathrm{tr}\left(\mathrm{ad}^2(v)\right) + O\left(\Vert v \Vert^4\right)\right)dv
$$
where $G(v\hspace{0.03cm};\sigma) = (z(\sigma))^{-1}e^{-\Vert v\Vert^2/2\sigma^2}$ is a normal probability density. 
If $B:\mathfrak{g} \times \mathfrak{g} \rightarrow \mathbb{R}$ denotes the Killing form, 
it now follows that~\cite{besse} (Chapter 7, Page 194)
$$
Z(\sigma) = z(\sigma)\,\int_{T_{x_o}M}G(v\hspace{0.03cm};\sigma)\left(1 + \frac{1}{12}B(v,v)+ O\left(\Vert v \Vert^4\right)\right)dv
$$
Finally, using basic properties of normal distributions
\begin{equation} \label{eq:prooftaylor3}
Z(\sigma) = z(\sigma)\left(1 + \frac{\sigma^2}{12}\mathrm{tr}(B_{\mathfrak{p}})+ O\left(\sigma^4\right)\right)
\end{equation}
where $B_{\mathfrak{p}}$ is the restriction of $B$ to $\mathfrak{p} \times \mathfrak{p}$ and its trace is taken with respect to the Riemannian metric of $M$. It~is~clear that (\ref{eq:prooftaylor3}) concludes the proof of (\ref{eq:Ztaylor}).\\[0.1cm]
\noindent \textbf{Remark\,:} it is interesting to compare the expansion (\ref{eq:prooftaylor3}) to a similar expansion that follows from (\ref{eq:typeivZ_bis}) and (\ref{eq:typeivz_bis}). These two identities yield
\begin{equation} \label{eq:prooftaylor4}
Z(\sigma) = z(\sigma)\left(1 + \sigma^2(\rho,\rho)+ O\left(\sigma^4\right)\right)
\end{equation}
Then, (\ref{eq:prooftaylor3}) and (\ref{eq:prooftaylor4}) imply that $\mathrm{tr}(B_{\mathfrak{p}}) = 12\,(\rho,\rho)$. One may always choose the Riemannian metric equal to $B_\mathfrak{p}\hspace{0.03cm}$. With this choice, $\mathrm{tr}(B_{\mathfrak{p}}) = d$ (the dimension of $M$), and 
\begin{equation} \label{eq:freudenthal}
   \dim\,M = 12\hspace{0.03cm}(\rho,\rho)
\end{equation}
In a slightly different notation, this is the well-known strange formula of Freudenthal and de Vries~\cite{freudenthal} (specifically, there is a factor $12$ instead of the usual $24$, because $\mathfrak{g}$ is here the realification of a simple complex Lie algebra, not~that simple complex Lie algebra itself). 

\subsection{Proof of Proposition \ref{prop:typeiv}}
Clearly, (\ref{eq:typeivZ}) is an immediate consequence of (\ref{eq:typeiv}). The proof of (\ref{eq:typeiv}) relies on the following two identities. First~\cite{helgasonbis}
(Chapter II, Page 268)
\begin{equation} \label{eq:iv_vol}
\prod_{\lambda \in \Delta_+} 2\sinh \lambda(a) = \sum_{w \in W}(\det w)e^{\rho(wa)}\;\;\;\;\;\;\;
\end{equation}
where the sum is over the elements $w$ of the Weyl group $W$, which are orthogonal transformations of $\mathfrak{a}$. Second~\cite{helgasonbis} (Chapter IV, Page 432)
\begin{equation} \label{eq:iv_spherical}
  \Phi_\tau(\mathrm{Exp}(a)) = \frac{\varpi(\rho)}{\varpi(\tau)}\frac{\sum_{w \in W} (\det w) e^{\tau(wa)}}
{\sum_{w \in W} (\det w) e^{\rho(wa)}}
\end{equation}
where $\mathrm{Exp}$ denotes the Riemannian exponential map. The idea is to replace these identities into the integral formula
\begin{align}
\nonumber  Z(\sigma,\tau) & = \int_M\,e^{-d^{\hspace{0.02cm} 2}(x,x_o)/2\sigma^2}\hspace{0.02cm}\Phi_\tau(x)\hspace{0.02cm}\mathrm{vol}(dx) \\[0.12cm]
\label{eq:proof22}                       & = \frac{\omega_M}{|W|}\int_{\mathfrak{a}} e^{-\Vert a\Vert^2/2\sigma^2}\hspace{0.02cm}\Phi_\tau(\mathrm{Exp}(a))\prod_{\lambda \in \Delta_+} \sinh^2\lambda(a)\hspace{0.02cm}da  
\end{align}
which is (just like (\ref{eq:Zcartan})) a consequence of the integral formula for the Cartan decomposition~\cite{helgasonbis} (Page 186). Note~that, since $M$ is a type IV symmetric space, all multiplicities $m_\lambda$ are equal to $2$, so all $\sinh \lambda(a)$ appear with a power $2$. \hfill\linebreak
Replacing (\ref{eq:iv_vol}) and (\ref{eq:iv_spherical}) into (\ref{eq:proof22}), it immediately follows that
$$
\nonumber  Z(\sigma,\tau) = \frac{\omega_M}{4^{|\Delta_+|}}\times \frac{\varpi(\rho)}{\varpi(\tau)}\times \frac{1}{|W|}
\sum_{w \in W}\sum_{w^\prime \in W} (\det ww^\prime)
\int_{\mathfrak{a}} e^{-\Vert a\Vert^2/2\sigma^2}\hspace{0.02cm}e^{(w^\prime\tau+w\rho)(a)}\hspace{0.02cm}da
$$
By the formula for the moment generating function of a normal distribution, this is 
\begin{equation} \label{eq:proof23}
Z(\sigma,\tau) = \frac{\omega_M}{4^{|\Delta_+|}}\times \frac{\varpi(\rho)}{\varpi(\tau)}\times (2\pi\sigma^2)^{r/2}\times \frac{1}{|W|}
\sum_{w \in W}\sum_{w^\prime \in W} (\det ww^\prime)
\exp\left[\frac{\sigma^2}{2}\Vert w^\prime\tau+w\rho\Vert^2\right]
\end{equation}
which can be further simplified by noting that
$$
\Vert w^\prime\tau+w\rho\Vert^2 = \Vert w^\prime\tau\Vert^2 + \Vert w \rho\Vert^2 + 2(w\rho,w^\prime\tau) = 
                                                   \Vert \tau\Vert^2 + \Vert \rho\Vert^2 + 2(w\rho,w^\prime\tau)
$$
since $w$ and $w^\prime$ are orthogonal. Substituting into (\ref{eq:proof23}), 
\begin{equation} \label{eq:proof24}
Z(\sigma,\tau) = \frac{\omega_M}{4^{|\Delta_+|}}\times \frac{\varpi(\rho)}{\varpi(\tau)}\times (2\pi\sigma^2)^{r/2} \times e^{\frac{\sigma^2}{2}\left(\Vert\tau\Vert^2+\Vert\rho\Vert^2\right)}\times \frac{1}{|W|}
\sum_{w \in W}\sum_{w^\prime \in W} (\det ww^\prime)
e^{\sigma^2(w\rho,w^\prime\tau)}
\end{equation}
Now, in order to prove (\ref{eq:typeiv}), it remains to show
\begin{equation} \label{eq:proof25}
\frac{1}{|W|}
\sum_{w \in W}\sum_{w^\prime \in W} (\det ww^\prime)
e^{\sigma^2(w\rho,w^\prime\tau)} = 2^{|\Delta_+|}\prod_{\lambda \in \Delta_+} \sinh\left(\sigma^2(\lambda,\tau)\right)
\end{equation}
However, this can be done by writing
$$
\sum_{w \in W}\sum_{w^\prime \in W} (\det ww^\prime)e^{\sigma^2(w\rho,w^\prime\tau)} = \sum_{w \in W}\sum_{w^\prime \in W} (\det w^{-1}w^\prime)e^{\sigma^2(\rho,w^{-1}w^\prime\tau)} = 
\sum_{w \in W}\sum_{w^\prime \in W} (\det w^\prime)e^{\sigma^2(\rho,w^\prime\tau)}
$$
where the first equality follows because $w$ is orthogonal and the second one by renaming $w^{-1}w^\prime$ to $w^\prime$. Indeed, (\ref{eq:proof25}) is then obtained by applying (\ref{eq:iv_vol}) to the last sum over $w^\prime$. Finally, (\ref{eq:proof24}) and (\ref{eq:proof25}) yield
$$
Z(\sigma,\tau) = \frac{\omega_M}{2^{|\Delta_+|}}\times \frac{\varpi(\rho)}{\varpi(\tau)}\times (2\pi\sigma^2)^{r/2} \times e^{\frac{\sigma^2}{2}\left(\Vert\tau\Vert^2+\Vert\rho\Vert^2\right)}\times \prod_{\lambda \in \Delta_+} \sinh\left(\sigma^2(\lambda,\tau)\right)
$$
which is essentially the same as (\ref{eq:typeiv}).

\subsection{Proof of Proposition \ref{prop:free}} 
Dividing (\ref{eq:typeivZ_bis}) by (\ref{eq:typeivz_bis}) and taking logarithms, one has the identity
\begin{equation} \label{eq:asymptotics_start}
  \frac{1}{r^2}\log\!\left(\frac{Z(\sigma)}{z(\sigma)}\right)   = \frac{t}{r^3}\hspace{0.03cm}(\rho,\rho) \,+\, \frac{1}{r^2}\sum_{\lambda \in \Delta_+} \log\!\left(\mathrm{sch}\left((t/r)(\lambda,\rho)\right)\right) 
\end{equation}
The left-hand side will now be denoted $F_r(t)$ and considered for each one of the classical reduced root systems. \\[0.1cm]
\noindent $A_{r}$ root system\,: here,
$\mathfrak{a}$ is the space of $(r+1) \times (r+1)$ real diagonal matrices $a = \mathrm{diag}(a_j)$ with $\mathrm{tr}(a) = 0$.\hfill\linebreak The positive roots are $\lambda_{ij}(a) = a_i - a_j$ for $i < j$ so that $\rho(a) = 2 \sum_j (r+1-j)\hspace{0.02cm} a_j\hspace{0.03cm}$. If~the Riemannian metric is given by $\Vert a \Vert^2 = 4\hspace{0.04cm}\mathrm{tr}(a^2)$, then the roots are normalised as in the proposition. Moreover, (\ref{eq:asymptotics_start}) implies 
$$
F_r(t) =  \frac{t}{r^3} \sum^r_{j=1} j^2  + \frac{1}{r^2} \sum_{i<j} \log\!\left[\mathrm{sch}
\left(\frac{t}{2}\left(\frac{j-i}{r}\right)\right)\right]
$$
The sum in the first term is elementary, while the second term is one half of a double Riemann sum. 
Specifically, because $\log(\mathrm{sch}(z))$ is a smooth function of real $z$, equal to $0$ when $z = 0$, 
\begin{equation} \label{eq:proof_asymptotics_A}
  F_r(t) = \frac{t}{3} + \frac{1}{2}\,\int^1_0\int^1_0 \log\left(\mathrm{sch}\left((t/2)(y-x)\right)\right)dx\hspace{0.02cm}dy + O\left(\frac{1}{r}\right)
\end{equation}
To obtain (\ref{eq:asymptotics_A}), note that for~any even continuous function $\phi$, 
$$
\int^1_0 \int^1_0 \phi(y-x)\hspace{0.02cm}dxdy = 2\,\int^1_0(1-x)\phi(x)\hspace{0.02cm}dx
$$
and apply this identity to (\ref{eq:proof_asymptotics_A}) with $\phi(x) = \log\left(\mathrm{sch}\left(tx/2\right)\right)$. \\[0.1cm]
\noindent $B_r$ root system\,: $\mathfrak{a}$ is the space of $r \times r$ real diagonal matrices $a = \mathrm{diag}(a_j)$ and the positive roots are $\lambda^-_{ij}(a) = a_i - a_j$ and $\lambda^+_{ij}(a) = a_i+a_j$ for $i < j$, in addition to $\lambda_j(a) = a_j$ for $j=1,\ldots,r$, so that $\rho(a) = \sum_j(2r-2j+1)a_j$. If~the Riemannian metric is given by $\Vert a \Vert^2 = 4\hspace{0.04cm}\mathrm{tr}(a^2)$, then the roots are normalised as required in the proposition. Moreover, (\ref{eq:asymptotics_start}) becomes    
\begin{align} 
\nonumber F_r(t)  =  \frac{t}{4r^3} \sum^r_{j=1} (2j-1)^2  &+ \frac{1}{r^2} \sum_{i<j} \log\!\left[\mathrm{sch}
\left(\frac{t}{2}\left(\frac{j-i}{r}\right)\right)\right] \\[0.1cm]
\nonumber  & + \frac{1}{r^2} \sum_{i<j}
\log\!\left[\mathrm{sch}
\left(\frac{t}{2}\left(2 - \frac{j+i}{r}+ \frac{1}{r}\right)\right)\right] + \frac{1}{r^2} \sum^r_{j=1}
\log\!\left[\mathrm{sch}
\left(\frac{t}{2}\left(1- \frac{j}{r}+\frac{1}{2r}\right)\right)\right]
\end{align}
This can be simplified just as in the previous case of $A_r$ root systems (\textit{i.e.} by identifying double Riemann sums).
In this way, and after noting that the last term on the right-hand side is $1/r^2$ times a simple (not a double) sum,
\begin{align} \label{eq:proof_asymptotics_B}
\nonumber  F_r(t) = \frac{t}{3} & + \frac{1}{2}\,\int^1_0\int^1_0 \log\left(\mathrm{sch}\left((t/2)(y-x)\right)\right)dx\hspace{0.02cm}dy \\[0.1cm]
& + \frac{1}{2}\,\int^1_0\int^1_0 \log\left(\mathrm{sch}\left((t/2)(y+x)\right)\right)dx\hspace{0.02cm}dy
+ O\left(\frac{1}{r}\right)
\end{align}
To obtain (\ref{eq:asymptotics_B}), note that for any even continuous function $\phi$,
$$
\int^1_0\int^1_0\phi(y-x)\hspace{0.02cm}dxdy + \int^1_0\int^1_0\phi(y+x)\hspace{0.02cm}dxdy = \int^2_0(2-x)\phi(x)\hspace{0.02cm}dx  
$$
and apply this to (\ref{eq:proof_asymptotics_B}) with $\phi(x) = \log\left(\mathrm{sch}\left(tx/2\right)\right)$. \\
For $C_r$ and $D_r$ root systems, the proof is mostly identical. It is enough to replace into (\ref{eq:asymptotics_start}) the following information. \\[0.1cm]
\noindent $C_{r}$ root system\,:
$\mathfrak{a}$ is the space of $r \times r$ real diagonal matrices $a = \mathrm{diag}(a_j)$ and the positive roots are 
$\lambda^-_{ij}(a) = a_i - a_j$ and $\lambda^+_{ij}(a) = a_i+a_j$ for $i < j$, in addition to $\lambda_j(a) = 2a_j$ for $j=1,\ldots,r$, so that $\rho(a) = 2 \sum_j (r+1-j)\hspace{0.02cm} a_j\hspace{0.03cm}$. 
The Riemannian metric should be fixed as before, $\Vert a \Vert^2 = 4\hspace{0.04cm}\mathrm{tr}(a^2)$. \\[0.1cm]
\noindent $D_{r}$ root system\,: $\mathfrak{a}$ is the space of $r \times r$ real diagonal matrices $a = \mathrm{diag}(a_j)$ and the positive roots are 
$\lambda^-_{ij}(a) = a_i - a_j$ and $\lambda^+_{ij}(a) = a_i+a_j$ for $i < j$. Then, $\rho(a) = 2 \sum_j (r-j)\hspace{0.02cm} a_j\hspace{0.03cm}$. 
The Riemannian metric is once more $\Vert a \Vert^2 = 4\hspace{0.04cm}\mathrm{tr}(a^2)$.

\subsection{Proof of Proposition \ref{prop:equilibrium}}
The proof follows closely in the footsteps of~\cite{deift} (Chapter 6). It is divided into four parts.

\subsubsection{Equilibrium measures}  
The aim here is to prove the existence of so-called equilibrium measures. These are probability measures on $\mathbb{R}$\,:
$\mu^*_{A}$, $\mu^*_{C}$ and $\mu^*_{BC}$ which achieve $E_{A} = I_{A}(\mu^*_{A})$, $E_{C} = I_{C}(\mu^*_{C})$ and $E_{BC} = I_{BC}(\mu^*_{BC})$. No claim will be made as to their uniqueness. 
\begin{proposition} \label{prop:equilibrium1}
  The equilibrium energies $E_A$, $E_C\hspace{0.01cm}$, and $E_{BC}$ are finite, and there exist probability measures 
 which achieve $E_{A} = I_{A}(\mu^*_{A})$, $E_{C} = I_{C}(\mu^*_{C})$ and $E_{BC} = I_{BC}(\mu^*_{BC})$. Moreover, any one of these probability measures is compactly supported. 
\end{proposition}
This proposition will be obtained by applying the following lemma.
\begin{lemma} \label{lem:equilibrium1}
  For each probability measure $\mu$ on $\mathbb{R}$ consider the functional
\begin{equation} \label{eq:generalFI}
    I(\mu) = \int\int k(x,y)\hspace{0.02cm}\mu(dx)\mu(dy)
\end{equation}
and let $E = \inf_\mu I(\mu)$. Assume that $k:\mathbb{R} \times \mathbb{R} \rightarrow \mathbb{R} \cup \lbrace +\infty \rbrace$ satisfies\,:
\begin{enumerate}
    \item[(a)]  There exists a continuous function $\psi : \mathbb{R} \rightarrow \mathbb{R}$ such that
$$
k(x,y) \geq \frac{1}{2}\psi(x) + \frac{1}{2}\psi(y)
$$
Moreover, $\psi(x) > c_\psi > -\infty$ for any real $x$, and $\psi(x) \rightarrow +\infty$ as $|x| \rightarrow \infty$. 
   \item[(b)] There exists a probability measure $\mu$ on $\mathbb{R}$ such that $I(\mu) < +\infty$. 
   \item[(c)] The function $\min\lbrace k,H\rbrace : \mathbb{R} \times \mathbb{R} \rightarrow \mathbb{R}$ (pointwise minimum) is continuous for any $H > 0$. 
\end{enumerate}
Then, $E$ is finite ($-\infty < E < + \infty$) and there exists at least one probability measure $\mu^*$ on $\mathbb{R}$ such that $E = I(\mu^*)$. Moreover, any such $\mu^*$ is compactly supported.  
\end{lemma}
To obtain Proposition \ref{prop:equilibrium1}, it is enough to cast (\ref{eq:functionala})-(\ref{eq:functionalbc}) in the form (\ref{eq:generalFI}) and then verify the conditions in the lemma. 
Each one of (\ref{eq:functionala})-(\ref{eq:functionalbc}) is indeed of the form (\ref{eq:generalFI}), with $k$ given by
\begin{equation} \label{eq:specialk}
k(x,y) = \frac{1}{2}W(x) + \frac{1}{2}W(y) - \frac{\beta}{2} L(x,y)
\end{equation}
for the following functions\,:
\begin{align}
\label{eq:specialka} \text{for (\ref{eq:functionala})}   && W(x) = \frac{c}{2t}x^2 && L(x,y) = \ell(x-y) \\[0.1cm]
\label{eq:specialkc} \text{for (\ref{eq:functionalc})}   && W(x) = \frac{c}{2t}x^2 && L(x,y) = \ell(x-y) + \ell(x+y) \\[0.1cm]
 \label{eq:specialkbc} \text{for (\ref{eq:functionalbc})} && W(x) = \frac{c}{2t}\hspace{0.02cm}x^2 - \beta(\delta-1)\hspace{0.02cm}\ell(x) && L(x,y) = \ell(x-y) + \ell(x+y) 
\end{align}
In each one of these three cases, Condition (a) is verified with 
$$
\psi(x) = \frac{c}{2t}x^2 - \beta(\delta+1)\ell(2|x|+1)
$$
which is indeed continuous and goes to $+\infty$ as $|x|\rightarrow \infty$. In particular, there exists $c_\psi > -\infty$ such that $\psi(x) > c_\psi$ for all $x \in \mathbb{R}$. Condition (b) is verified for $\mu$ a standard Gaussian distribution, and Condition (c) follows from the elementary fact that $\ell(x) < +\infty$ for any $x \in \mathbb{R}$.   

Thus, in order to prove Proposition \ref{prop:equilibrium1} it remains to prove Lemma \ref{lem:equilibrium1}. \\[0.1cm]
\textbf{Proof of Lemma \ref{lem:equilibrium1}\,:} Condition (a) implies that $k(x,y) > c_\psi$ for all $(x,y) \in \mathbb{R} \times \mathbb{R}$, so $I(\mu) > c_\psi$ for any probability measure $\mu$, and $E > c_\psi > -\infty$. On the other hand, Condition (b) implies that $E < + \infty$. Therefore, $E$ is finite. \hfill\linebreak
\indent To prove the existence of $\mu^*$ such that $E = I(\mu^*)$, let $(\mu_n;n=1,2,\ldots)$ be a sequence of probability measures such that $I(\mu_n) \rightarrow E$. It is now shown this sequence is tight. In particular~\cite{billingsley}, this implies it has a subsequence which converges weakly to some probability measure $\mu^*$.

To show that tightness holds, note that for any $\varepsilon > 0$ one has $I(\mu_n) < E + \varepsilon$ for all large enough $n$. Using Condition (a), it follows that
$$
E + \varepsilon > \int\int k(x,y)\hspace{0.02cm}\mu_n(dx)\mu_n(dy) \geq \frac{1}{2}\int \psi(x)\hspace{0.02cm}\mu_n(dx) + 
\frac{1}{2}\int \psi(y)\hspace{0.02cm}\mu_n(dy)
$$
In other words,
$$
E + \varepsilon > \int \psi(x)\hspace{0.02cm}\mu_n(dx)
$$
Since the left-hand side does not depend on $n$, and $\psi(x) \rightarrow +\infty$ as $|x| \rightarrow \infty$, the sequence $(\mu_n)$ is indeed tight. 

Next, it is shown that $E = I(\mu^*)$. To do so, let $(\mu_{n(k)})$ be a subseuqnce of $(\mu_n)$ which converges weakly to $\mu^*$. This implies the product measures $\mu_{n(k)}\otimes \mu_{n(k)}$ converge weakly to $\mu^* \otimes \mu^*$~\cite{billingsley} (see Theorem 2.8, Page 23). By~Conditions (a) and (c), $\min\lbrace k,H\rbrace$ is a bounded continuous function, so weak convergence implies the limit
$$
\int \int \min\lbrace k(x,y),H\rbrace\hspace{0.02cm} \mu_{n(k)}(dx)\mu_{n(k)}(dy) \rightarrow \int \int \min\lbrace k(x,y),H\rbrace\hspace{0.02cm} \mu^*(dx)\mu^*(dy) 
$$
By recalling the definition (\ref{eq:generalFI}) of $I(\mu_{n(k)})$, it becomes clear
$$
\liminf I(\mu_{n(k)}) \geq \int \int \min\lbrace k(x,y),H\rbrace\hspace{0.02cm} \mu^*(dx)\mu^*(dy)
$$
Then, since $H > 0$ is arbitrary, using the monotone convergence theorem, one has $\liminf I(\mu_{n(k)}) \geq I(\mu^*)$.
However, $I(\mu_{n(k)}) \rightarrow E$ and therefore $E \geq I(\mu^*)$, which shows $E = I(\mu^*)$. 

It only remains to show that $\mu^*$ is compactly supported. Let $D \subset \mathbb{R}$ and assume that $\mu^*(D) > 0$. Then, define 
$$
\tilde{\mu}_{\varepsilon} = (1+\varepsilon \mu^*(D))^{-1}\hspace{0.02cm}(\mu^* + \varepsilon \mu^*|_D)
$$
where $\mu^*|_D$ is the restriction of $\mu^*$ to $D$. Note that $\tilde{\mu}_\varepsilon$ is a probability measure if 
$|\varepsilon| < 1$. Moreover,
$$
I(\tilde{\mu}_{\varepsilon}) = (1+\varepsilon \mu^*(D))^{-2}\int\int k(x,y)\hspace{0.02cm}(\mu^* + \varepsilon \mu^*|_D)(dx) 
(\mu^* + \varepsilon \mu^*|_D)(dy) 
$$
reaches a minimum at $\varepsilon = 0$. In particular, $d/d\varepsilon|_{\varepsilon = 0}\hspace{0.04cm}I(\tilde{\mu}_{\varepsilon}) = 0$, which means (by a straightforward calculation)
$$
\int \int (k(x,y) - I(\mu^*))\hspace{0.02cm}\mu^*|_D(dx)\mu^*(dy) = 0
$$
Applying Condition (a) to this last identity,
$$
0 \geq \int \int (\psi(x)/2 + \psi(y)/2 - I(\mu^*))\hspace{0.02cm}\mu^*|_D(dx)\mu^*(dy)
$$
but, after integrating with respect to $\mu^*(dy)$, this is the same as
$$
0 \geq \int \left[\psi(x) + \int \psi(y)\mu^*(dy) - 2I(\mu^*)\right]\mu^*|_D(dx) 
$$
Since $\psi(x) \rightarrow +\infty$ as $|x| \rightarrow \infty$, there exists $M > 0$ such that the integrand is $\geq 1$ for $|x| > M$. If $D$ is contained in the complement of $[-M,M]$, one must therefore have $\mu^*(D) = 0$, in contradiction with the assumption $\mu^*(D) > 0$.

\subsubsection{$r$-point energy} Let $k$ be a function as in Lemma \ref{lem:equilibrium1}, and define
\begin{equation} \label{eq:KK}
  K_r(a) = \sum_{i \neq j} k(a_i,a_j) \hspace{0.5cm} \text{for $a \in \mathbb{R}^r$}
\end{equation}
where $r$ is a positive integer (in the sum, $1 \leq i,j \leq r$). Then, let 
\begin{equation} \label{eq:transdiameter}
  E_r = \inf_{a \in \mathbb{R}^r}\hspace{0.03cm}\frac{1}{r(r-1)}\hspace{0.02cm}K_r(a)
\end{equation}
This will be denoted by $E_r(A)$, $E_r(C)$, or $E_r(BC)$, when $k$ is given by (\ref{eq:specialk}) and (\ref{eq:specialka}), (\ref{eq:specialkc}), or (\ref{eq:specialkbc}), respectively. These $E_r(A)$, $E_r(C)$, and $E_r(BC)$ will be called $r$-point energies. The aim here will be to prove the proposition.
\begin{proposition} \label{prop:fekete}
  The $r$-point energies $E_r(A)$, $E_r(C)$, and $E_r(BC)$ are finite, and converge to the corresponding equilibrium energies $E_A$, $E_C$, and $E_{BC}$ as $r \rightarrow \infty$ (that is, $E_r(A) \rightarrow E_A$, \textit{etc.}). 
\end{proposition}
This will be obtained by applying the following lemma.
\begin{lemma} \label{lem:fekete}
  Let $k: \mathbb{R} \times \mathbb{R} \rightarrow \mathbb{R} \cup \lbrace + \infty \rbrace$ satisfy all the conditions of Lemma \ref{lem:equilibrium1}, and define $E_r$ as in (\ref{eq:transdiameter}). If~$k(x,x) = +\infty$ for all $x \in \mathbb{R}$, then $E_r$ is finite and $E_r \rightarrow E$ as $r \rightarrow \infty$.  
\end{lemma}
Admitting this lemma, Proposition \ref{prop:fekete} becomes immediate. Indeed, when $k$ is given by (\ref{eq:specialk}) and (\ref{eq:specialka}), (\ref{eq:specialkc}), or (\ref{eq:specialkbc}), it is clear that $k(x,x) = +\infty$ for any $x \in \mathbb{R}$. \\[0.1cm]
\textbf{Proof of Lemma \ref{lem:fekete}\,:} Condition (a) from Lemma \ref{lem:equilibrium1} implies $K_r(a) > r(r-1)c_\psi$ so $E_r > c_\psi\hspace{0.03cm}$. On the other hand, for any probability measure $\mu$, it holds that
\begin{equation} \label{eq:multipleK}
\int\ldots\int K_r(a_{\scriptscriptstyle 1},\ldots,a_r)\hspace{0.02cm}\mu(da_{\scriptscriptstyle 1})\ldots\mu(da_r) = r(r-1) \int\int k(x,y)\hspace{0.02cm}\mu(dx)\mu(dy)
\end{equation}
where $K_r(a_{\scriptscriptstyle 1},\ldots,a_r)$ denotes $K_r(a)$ when $a = (a_{\scriptscriptstyle 1},\ldots,a_r)$. Indeed, (\ref{eq:multipleK}) follows by replacing (\ref{eq:KK}) under the multiple integral and then integrating over each pair of coordinates $(a_i,a_j)$, since there are $r(r-1)$ such pairs. Note~that the right-hand side of (\ref{eq:multipleK}) is equal to $r(r-1)\hspace{0.02cm}I(\mu)$, by (\ref{eq:generalFI}). Then, taking $\mu = \mu^*$ where $E = I(\mu^*)$ (as~in Lemma \ref{lem:equilibrium1}), it follows that $\inf_a K_r(a) \leq r(r-1)E$. In other words, $E_r \leq E$, so that $E_r$ is indeed finite. \hfill\linebreak
\indent Conditions (a) and (c) from Lemma \ref{lem:equilibrium1} ensure that $K_r$ is bounded below, that $K_r(a) \rightarrow +\infty$ when $|a| \rightarrow \infty$, and that $K_r$ is lower semicontinuous (by Condition (c), it is a pointwise supremum of continuous functions). Therefore, the infimum in (\ref{eq:transdiameter}) is achieved\,: there exists $a^* = (a^*_{\scriptscriptstyle 1},\ldots,a^*_{\scriptscriptstyle r})$ such that $E_r = (r(r-1))^{-1}\hspace{0.02cm}K_r(a^*)$. 

To prove that $E_r \rightarrow E$, consider probability measures $\nu_r = \frac{1}{r}\sum^r_{i=1} \delta_{a_i}$ ($\delta_{a_i}$ is the Dirac measure at $a_i \in \mathbb{R}$). These~probability measures are tight. Indeed, by (\ref{eq:transdiameter}) and Condition (a) from Lemma \ref{lem:equilibrium1},
$$
E_r \geq \frac{1}{r(r-1)}\sum_{i\neq j} \frac{1}{2}\psi(a^*_i) + \frac{1}{2} \psi(a^*_j) = \frac{1}{r(r-1)}\times (r-1) \sum^r_{i=1} \psi(a^*_i)
$$
and (since $E \geq E_r$) this means that 
$$
E \geq \int\psi(x)\nu_n(dx)
$$
Tightness of $(\nu_r)$ therefore follows since $\psi(x) \rightarrow \infty$ as $|x| \rightarrow \infty$. Then, by reducing to a subsequence if necessary, it is possible to assume that $(\nu_r)$ converges weakly to some probability measure $\nu$. 
Now, for any $H > 0$, note that
$$
\int\int \min\lbrace k(x,y),H\rbrace\hspace{0.02cm}\nu_r(dx) \nu_r(dy) = \frac{1}{r^2} \sum^r_{i,j=1} \min\lbrace k(a^*_i,a^*_j),H\rbrace = \frac{1}{r^2} \sum_{i \neq j} \min\lbrace k(a^*_i,a^*_j),H\rbrace + \frac{H}{r}
$$
where the second equality follows from $k(a^*_i,a^*_i) = +\infty$ for $i = 1,\ldots,r$. Multiplying both sides by $r^2/(r(r-1))$,
letting $r \rightarrow \infty$, and using the fact that $\min\lbrace k,H\rbrace$ is a bounded continuous function (Conditions (a) and (c) from~Lemma \ref{lem:equilibrium1}) and $\nu_r$ converge weakly to  $\nu$, this yields 
$$
\int\int \min\lbrace k(x,y),H\rbrace\hspace{0.02cm}\nu(dx) \nu(dy) \leq \liminf E_r
$$
However, since $H$ is arbitrary, it now follows from (\ref{eq:generalFI}) that $I(\nu) \leq \liminf E_r\hspace{0.03cm}$. Finally, since $E = I(\mu^*) \leq I(\nu)$ and $E_r \leq E$, it follows that $E \leq \liminf E_r \leq \limsup E_r \leq E$, and therefore $E_r \rightarrow E$. \\[0.1cm]
\noindent \textbf{Remark\,:} in fact, the above shows that $I(\nu) = E$, so that $\nu$ is an equilibrium measure. 
\subsubsection{Stability under smoothing} Let $k$, $I$ and $\mu^*$ be as in Lemma \ref{lem:equilibrium1}. Denote by $\mu^*_\varepsilon$ the smoothed version of $\mu^*$, which has the following density with respect to the Lebesgue measure on $\mathbb{R}$, for $\varepsilon > 0$,
\begin{equation} \label{eq:smoothmeasure}
  m^*_\varepsilon(x) = \frac{1}{2\varepsilon}\int^{x+\varepsilon}_{x-\varepsilon}\mu^*(dy)
\end{equation}
Specifically, $\mu^*_\varepsilon(dx) = m^*_\varepsilon(x)\hspace{0.02cm}dx$. The aim here will be to prove the following proposition. 
\begin{proposition} \label{prop:stability}
 In the setting of Proposition \ref{prop:equilibrium1}, let $I$ denote $I_A$, $I_C$ or $I_{BC}$ and let $\mu^*$ denote a corresponding equilibrium measure $\mu^*_A$, $\mu^*_C$ or $\mu^*_{BC}\hspace{0.03cm}$. Then, for the smoothed measure $\mu^*_\varepsilon$, it holds that $I(\mu^*_\varepsilon) \rightarrow I(\mu^*)$ as $\varepsilon \rightarrow 0$. 
\end{proposition}
Roughly, this proposition says that equilibrium is stable under smoothing. It will be proved using a general lemma.
\begin{lemma} \label{lem:stability}
   Let $k$, $I$ and $\mu^*$ be as in Lemma \ref{lem:equilibrium1}. Assume that $k$ satisfies all the conditions of Lemmas \ref{lem:equilibrium1} and \ref{lem:fekete}. In~addition, assume that $k$ can be written
\begin{equation} \label{eq:stablecondition}
  k(x,y) = h(x,y) - \Lambda(x,y)
\end{equation}
where $h:\mathbb{R} \times \mathbb{R} \rightarrow \mathbb{R}$ is continuous and $\Lambda(x,y)$ is a logarithmic part
\begin{equation} \label{eq:Lambda}
  \Lambda(x,y) = \sum^{P}_{p=1} a_p\hspace{0.02cm}\log|l_p(x,y)| 
\end{equation}
with $a_p > 0$ and $l_p:\mathbb{R} \times \mathbb{R} \rightarrow \mathbb{R}$ a linear function for each $p = 1, \ldots, P$. Then, for the smoothed measure $\mu^*_\varepsilon$ (whose density is given by (\ref{eq:smoothmeasure})), it holds that $I(\mu^*_\varepsilon) \rightarrow I(\mu^*)$ as $\varepsilon \rightarrow 0$.
\end{lemma}
Once this lemma has been proved, Proposition \ref{prop:stability} can be obtained by showing that Condition (\ref{eq:stablecondition}) is verified when $k$ is given by (\ref{eq:specialk}) and (\ref{eq:specialka}), (\ref{eq:specialkc}), or (\ref{eq:specialkbc}). This is indeed the case,  with the following functions $h$ and $\Lambda$, 
\begin{align*}
\text{for (\ref{eq:specialka})}   && h(x,y) = \frac{c}{4t}(x^2+y^2) - \log\left(\mathrm{sch}(x-y)\right) && \Lambda(x,y) = \log|x-y| \\[0.12cm]
\text{for (\ref{eq:specialkc})}   && h(x,y) = \frac{c}{4t}(x^2+y^2) - \log\left(\mathrm{sch}(x-y)\right) - \log\left(\mathrm{sch}(x+y)\right) && \Lambda(x,y) = \log|x-y|+\log|x+y| \\[0.12cm]
\text{for (\ref{eq:specialkbc})} && h(x,y) = \frac{c}{4t}(x^2+y^2) - \log\left(\mathrm{sch}(x-y)\right) - \log\left(\mathrm{sch}(x+y)\right)  && \Lambda(x,y) = \log|x-y|+\log|x+y| \\[0.1cm]
&& - \frac{\beta}{2}(\delta-1)(\log\left(\mathrm{sch}(x)\right)+\log\left(\mathrm{sch}(y)\right)) && 
+ \frac{\beta}{2}(\delta-1)(\log|x|+\log|y|)
\end{align*}
where $\mathrm{sch}(z) = \sinh(z)/z$. The key point here is that $\log(\mathrm{sch}(z))$ is a continuous function of the real variable $z$. It now remains to prove Lemma \ref{lem:stability}. \\[0.1cm]
\textbf{Proof of Lemma \ref{lem:stability}\,:} The proof is broken down into several steps. \\[0.1cm]
\noindent \textbf{step 1\,:} recall that $I(\mu^*) = E$ is finite. By (\ref{eq:generalFI}) and (\ref{eq:stablecondition}), this is equal to 
\begin{equation} \label{eq:kmu*_11}
\int\int k(x,y)\hspace{0.02cm}\mu^*(dx)\mu^*(dy) = 
\int\int h(x,y)\hspace{0.02cm}\mu^*(dx)\mu^*(dy) -
\int\int \Lambda(x,y)\hspace{0.02cm}\mu^*(dx)\mu^*(dy)
\end{equation}
Lemma \ref{lem:equilibrium1} states that $\mu^*$ is compactly supported. Since $h$ is continuous, the first integral on the right-hand side converges. It follows that the second integral also converges. In other words, $-\Lambda$ is integrable with respect to the prodcut measure $\mu^*\otimes\mu^*$. \\[0.1cm]
\noindent \textbf{step 2\,:} since $k(x,x) = +\infty$ for all $x \in \mathbb{R}$, and since $k$ is integrable with respect to $\mu^*\otimes\mu^*$, the measure $\mu^*$ has no atoms. Therefore, the density $m^*_\varepsilon$ in (\ref{eq:smoothmeasure}) is continuous. Moreover, since $\mu^*$ is compactly supported, this density $m^*_\varepsilon$ is also compactly supported. This implies that $-\Lambda$ (as obtained from (\ref{eq:Lambda})) is integrable with respect to 
$\mu^*_\varepsilon\otimes\mu^*_\varepsilon\hspace{0.03cm}$.\hfill\linebreak By continuity of $h$, this shows that $I(\mu^*_\varepsilon)$ is finite.  \\[0.1cm]
\noindent \textbf{step 3\,:} consider then the expression of $I(\mu^*_\varepsilon)$. From (\ref{eq:generalFI}) and (\ref{eq:smoothmeasure}), this is equal to
$$
\int \int k(x,y)\hspace{0.02cm}m^*_\varepsilon(x)\hspace{0.02cm}m^*_\varepsilon(y)\hspace{0.02cm}dxdy = \frac{1}{4\varepsilon^2}\int\int k(x,y)\left(\int\int \mathds{1}_{[-\varepsilon,\varepsilon]}(u-x)\mathds{1}_{[-\varepsilon,\varepsilon]}(v-y)\mu^*(du)\mu^*(dv)\right)dx dy
$$
By changing the order of integration,
\begin{equation} \label{eq:littleboss}
I(\mu^*_\varepsilon) = \int \int \left(\frac{1}{4\varepsilon^2}\int^\varepsilon_{-\varepsilon}\int^\varepsilon_{-\varepsilon}
k(u-x,v-y)dxdy
 \right)\mu^*(du)\mu^*(dv) 
\end{equation}
The function inside the parenthesis will be denoted $k_\varepsilon(u,v)$. From (\ref{eq:stablecondition}), $k = h - \Lambda$, while $k_\varepsilon = h_\varepsilon - \Lambda_\varepsilon$ with $h_\varepsilon$~and~$\Lambda_\varepsilon$  obtained from $h$ and $\Lambda$ using the same ``smoothing" procedure. \\[0.1cm]
\noindent \textbf{step 4\,:} to prove that $I(\mu^*_\varepsilon) \rightarrow I(\mu^*)$, one has to prove 
\begin{equation} \label{eq:finalboss}
   \int\int \left[k_\varepsilon(u,v) - k(u,v)\right]\mu^*(du)\mu^*(dv) \rightarrow 0
\end{equation}
Of course, $k_\varepsilon - k = (h_\varepsilon - h) - (\Lambda_\varepsilon - \Lambda)$. Since $h$ is continuous, $h - h_\varepsilon \rightarrow 0$ uniformly on the support of $\mu^*\otimes \mu^*$. Therefore, (\ref{eq:finalboss}) is equivalent to 
\begin{equation} \label{eq:finalboss2}
   \int\int \left[\Lambda_\varepsilon(u,v) - \Lambda(u,v)\right]\mu^*(du)\mu^*(dv) \rightarrow 0
\end{equation}
\noindent \textbf{step 5\,:} before proving (\ref{eq:finalboss2}), note the following. Returning to (\ref{eq:Lambda}), one may write
$$
-\Lambda = \sum^P_{p=1} \mathds{1}_{\lbrace|l_p| \geq 1\rbrace}\times -a_p\log|l_p| + \sum^P_{p=1} \mathds{1}_{\lbrace|l_p| < 1\rbrace}\times -a_p\log|l_p|
$$
Here, each term in the first sum is a bounded continuous function on the compact support of $\mu^*\otimes \mu^*$. Integrability of $-\Lambda$ therefore implies integrability of the second sum (with respect to $\mu^*\otimes \mu^*$). However, each term in the second sum is positive. In conclusion, $-\log|l_p|$ is integrable with respect to $\mu^*\otimes \mu^*$, for each $p = 1, \ldots, P$. \\[0.1cm]
\noindent \textbf{step 6\,:} from (\ref{eq:Lambda}), note that $|\Lambda_\varepsilon - \Lambda|$ is less than the sum of the following $P$ terms,
\begin{equation} \label{eq:finalboss3}
\frac{1}{4\varepsilon^2}\int^\varepsilon_{-\varepsilon}\int^\varepsilon_{-\varepsilon} \left|\hspace{0.02cm}\log\left|l_p(u,v) - l_p(x,y)\right| - \log\left|l_p(u,v)\right|\hspace{0.02cm}\right|dxdy \hspace{1cm} (p = 1,\ldots,P)
\end{equation}
where the positive constants $a_p$ have been suppressed. After a linear transformation of the $(x,y)$-plane, $l_p(x,y) = x$, and
(\ref{eq:finalboss3}) becomes
$$
\frac{q}{2\varepsilon} \int^\varepsilon_{-\varepsilon} \left|\hspace{0.02cm}\log\left| 1 - (x/l)\right|\hspace{0.02cm}\right|dx  \leq q\log\left(1+(\varepsilon/l_p)\right) \\
$$
where $q > 0$ is a constant and $l_p$ denotes $l_p(u,v)$. After a simple manipulation of this last expression,
\begin{equation} \label{eq:finalboss4}
  \left|\Lambda_\varepsilon(u,v) - \Lambda(u,v) \right| \leq q^\prime \sum^P_{p=1} a_p\hspace{0.03cm}\left[\hspace{0.03cm}\log(l_p(u,v) + \varepsilon) - \log(l_p(u,v))\right]
\end{equation}
where $q^\prime > 0$ is a new constant. \\[0.1cm]
\noindent \textbf{step 7\,:} to prove (\ref{eq:finalboss2}), note from (\ref{eq:finalboss4}) that
$$
\int\int \left|\Lambda_\varepsilon(u,v) - \Lambda(u,v)\right|\mu^*(du)\mu^*(dv) \leq 
q^\prime \sum^P_{p=1} a_p\hspace{0.03cm}\int\int \left[\hspace{0.03cm}\log(l_p(u,v) + \varepsilon) - \log(l_p(u,v)) \right]\mu^*(du)\mu^*(dv)
$$
From step 5, $-\log|l_p|$ is integrable with respect to $\mu^*\otimes \mu^*$. Therefore, by the monotone convergence theorem, each integral on the right-hand side converges to $0$ as $\varepsilon \rightarrow 0$. Now, (\ref{eq:finalboss2}) follows immediately from the above inequality. 

\subsubsection{Conclusion of the proof} \label{subsec:conclusion_equilibrium}
The aim here is to finally prove Proposition \ref{prop:equilibrium}. To begin, recall the integrals $Z_\mathfrak{a}(\sigma)$ in (\ref{eq:IA})-(\ref{eq:IBC}) and the energy functionals $I_A$, $I_C$ and $I_{BC}$ in (\ref{eq:functionala})-(\ref{eq:functionalbc}). Each one of these functionals can be cast in the form (\ref{eq:generalFI}), by taking $k$ as in (\ref{eq:specialk}) and (\ref{eq:specialka})-(\ref{eq:specialkbc}). This~$k$~then satisfies all the conditions of Lemmas \ref{lem:equilibrium1}-\ref{lem:stability}. 

As of now, denote by $\tilde{Z}(r)$ the integral $Z_\mathfrak{a}(\sigma)$ with $\sigma^2 = t/r$. If $K_r$ is defined according to (\ref{eq:KK}), then $\tilde{Z}(r)$ can~be expressed in the following way.
\begin{itemize}
 \item[--] $A_{r}$ root system\,: if $Z_\mathfrak{a}(\sigma)$ is given by (\ref{eq:IA}),
\begin{equation} \label{eq:IAK}
\tilde{Z}(r) = \int_{\mathfrak{a}} \exp\left[-K_{r+1}(a)\hspace{0.02cm}\right]da
\end{equation}
where $k$ is given by (\ref{eq:specialk}) and (\ref{eq:specialka}).
 \item[--] $B_{r}$ root system\,: if $Z_\mathfrak{a}(\sigma)$ is given by (\ref{eq:IB}), with the additional assumption $\alpha = \beta(\delta-1)r$,
\begin{equation} \label{eq:IBK}
  \tilde{Z}(r) = \int_{\mathfrak{a}} \exp\left[-K_{r}(a)-\Psi_{r}(a)\right]da
\end{equation}
where $k$ is given by (\ref{eq:specialk}) and (\ref{eq:specialkbc}) and $\Psi_{r}(a) = \sum^{r}_{i=1} W(a_i)$.
\item[--] $C_{r}$ root system\,: if $Z_\mathfrak{a}(\sigma)$ is given by (\ref{eq:IC}), $\tilde{Z}(r)$ is given by (\ref{eq:IBK}),  with $k$ as in (\ref{eq:specialk}) and (\ref{eq:specialkc}), and with
$\Psi_r(a) = \sum^r_{i=1} W(a_i)- \gamma\hspace{0.03cm}\ell(2a_i)$.
 \item[--] $D_{r}$ root system\,: if $Z_\mathfrak{a}(\sigma)$ is given by (\ref{eq:ID}), $\tilde{Z}(r)$ is given by (\ref{eq:IBK}),  with $k$ as in (\ref{eq:specialk}) and (\ref{eq:specialkc}), and with $\Psi_r(a) = \sum^r_{i=1} W(a_i)$.
 \item[--] $(BC)_{r}$ root system\,: if $Z_\mathfrak{a}(\sigma)$ is given by (\ref{eq:IBC}), there are two subcases,
\begin{itemize}
   \item[(1)] with the additional assumption $\alpha = \beta(\delta-1)r$, $\tilde{Z}(r)$ is given by (\ref{eq:IBK}),  with $k$ as in (\ref{eq:specialk}) and (\ref{eq:specialkbc}), and $\Psi_r(a) = \sum^r_{i=1} W(a_i) - \gamma\hspace{0.03cm}\ell(2a_i)$. 
  \item[(2)] if $\alpha$ is independent of $r$ (last line of Table \ref{tab:classif}), $\tilde{Z}(r)$ is given by (\ref{eq:IBK}),  with $k$ as in (\ref{eq:specialk}) and (\ref{eq:specialkc}), and 
$\Psi_r(a) = \sum^r_{i=1} W(a_i) - \gamma\hspace{0.03cm}\ell(2a_i) - \alpha\hspace{0.03cm}\ell(a_i)$.
\end{itemize}           
\end{itemize}
Indeed, the above expressions of $\tilde{Z}(r)$ can be verified using the identity
\begin{equation} \label{eq:KKpotential}
K_r(a) = (r-1)\sum^r_{i=1} W(a_i) - \beta \sum_{i<j} L(a_i,a_j)
\end{equation}
which results from (\ref{eq:specialk}) and (\ref{eq:KK}). Replacing this identity into (\ref{eq:IAK}) and (\ref{eq:IBK}), it is straightforward to recover (\ref{eq:IA})-(\ref{eq:IBC}), with $\sigma$ rescaled according to $\sigma^2 = t/r$.

Proposition \ref{prop:equilibrium} will be obtained from the following general lemma.
\begin{lemma} \label{lem:laplace} Let $k$ satisfy all the conditions of Lemmas \ref{lem:equilibrium1}-\ref{lem:stability}, and define $K_r$ as in (\ref{eq:KK}). Also, let $w:\mathbb{R}\rightarrow \mathbb{R} \cup \lbrace +\infty\rbrace$ be a locally integrable function, and $\Psi_r(a) = \sum^r_{i=1} w(a_i)$ for any $a \in \mathbb{R}^r$ where $a = (a_{\scriptscriptstyle 1},\ldots,a_r)$. If $w$ satisfies
\begin{equation} \label{eq:laplace_condition}
    \int_{\mathbb{R}}\exp(-w(x))\hspace{0.03cm}dx < \infty
\end{equation}
then, the integrals
\begin{equation} \label{eq:laplace_Z}
  \mathcal{Z}(r) = \int_{\mathbb{R}^r}\exp\left[-K_{r}(a)-\Psi_{r}(a)\right]da
\end{equation}
satisfy the following limit,
\begin{equation} \label{eq:laplace}
\lim_{r\rightarrow \infty}\frac{1}{r^2}\log\!\left(\mathcal{Z}(r)\right) = -E
\end{equation}
where (the equilibrium energy) $E$ was defined in Lemma \ref{lem:equilibrium1}.
\end{lemma}
Interestingly, the limit (\ref{eq:laplace}) depends only on $k$ and ``forgets" the function $w$.

The application of Lemma \ref{lem:laplace} to Proposition \ref{prop:equilibrium} simply amounts to identifying $\tilde{Z}(r)$ from (\ref{eq:IAK}) or (\ref{eq:IBK}) with $\mathcal{Z}(r)$ from (\Ref{eq:laplace_Z}). In the latter case of (\ref{eq:IBK}), this identification is straightforward since $\mathfrak{a}$ is naturally identified with $\mathbb{R}^r$ (as~already explained after 
(\ref{eq:IBC})). In the case of (\ref{eq:IAK}), $\mathfrak{a}$ is the subspace (hyperplane) of $\mathbb{R}^{r+1}$ defined by the equation $\sum a_i = 0$ (sum over $i=1,\ldots,r+1$). However, also in this case, (\ref{eq:KKpotential}) implies
$$
K_{r+1}(a) = \frac{rc}{2t} \sum^{r+1}_{i=1} a^2_i - \beta \sum_{i<j} \ell(a_i-a_j)
$$
Therefore, letting $a_i = (r+1)^{-1/2}\hspace{0.03cm}x + \bar{a}_i$ where $x \in \mathbb{R}$ and $\bar{a} \in \mathfrak{a}$, the following decomposition can be obtained
$$
\int_{\mathbb{R}^{r+1}} \exp\left[-K_{r+1}(a)\hspace{0.02cm}\right]da = 
\left(\int_{\mathbb{R}}\exp\left[-\frac{rc}{2t}x^2\right]dx\right)
\left(\int_{\mathfrak{a}} \exp\left[-K_{r+1}(\bar{a})\hspace{0.02cm}\right]d\bar{a}\right)
$$
Now, this can be written $\mathcal{Z}(r+1) = \left(2\pi t/rc\right)^{\frac{1}{2}}\hspace{0.03cm}\tilde{Z}(r)$ (where $\mathcal{Z}(r)$ is given by (\ref{eq:laplace_Z}) with $\Psi_r(a) = 0$), which shows that $(1/r^2)\log\hspace{0.03cm}(\tilde{Z}(r))$ has the same limit as $(1/r^2)\log\hspace{0.03cm}(\mathcal{Z}(r))$ when $r \rightarrow \infty$. This limit is just the one in Lemma \ref{lem:laplace}. 

In conclusion, Proposition \ref{prop:equilibrium} follows because the limits (\ref{eq:lima})-(\ref{eq:limbc}) are of the form (\ref{eq:laplace}), where $E$ is equal to $E_A$,~$E_C$, or $E_{BC}\hspace{0.03cm}$,  according to whether $k$ is given by (\ref{eq:specialk}) and (\ref{eq:specialka}), (\ref{eq:specialkc}), or (\ref{eq:specialkbc}). \\[0.1cm]
\noindent \textbf{Proof of Lemma \ref{lem:laplace}\,:} First, it will be shown that
\begin{equation} \label{eq:proof_laplace_1}
 \limsup\, \frac{1}{r^2}\hspace{0.02cm}\log\!\left(\mathcal{Z}(r)\right) \leq -E
\end{equation}
To do so, note from (\ref{eq:laplace_Z}) that
$$
\mathcal{Z}(r) \leq \exp\left[-\inf_{a \in \mathbb{R}^r}\hspace{0.02cm}K_r(a)\right]\int_{\mathbb{R}^r}\exp\left[-\Psi_r(a)\right]da 
= \exp\left[-r(r-1)E_r\right]\left(\hspace{0.02cm}\int_{\mathbb{R}}\exp(-w(x))\hspace{0.03cm}dx\right)^{\!\!r}
$$
where $E_r$ was defined in (\ref{eq:transdiameter}). After taking logarithms,
$$
\frac{1}{r^2}\hspace{0.02cm}\log\!\left(\mathcal{Z}(r)\right) \leq \left(\frac{1}{r}-1\right)E_r + \frac{I_w}{r}
$$
where $I_w$ is the integral in (\ref{eq:laplace_condition}). Now, (\ref{eq:proof_laplace_1}) follows by letting $r \rightarrow \infty$ and noting that $E_r \rightarrow E$, by Lemma \ref{lem:fekete}. 
To~obtain the desired limit (\ref{eq:laplace}), it is then enough to show
\begin{equation} \label{eq:proof_laplace_2}
 \liminf\, \frac{1}{r^2}\hspace{0.02cm}\log\!\left(\mathcal{Z}(r)\right) \geq -E
\end{equation}
To do so, let $\mu^*$ be an equilibrium measure as in Lemma \ref{lem:equilibrium1}, and $\mu^*_\varepsilon$ the smoothed version of $\mu^*$, whose density $m^*_\varepsilon$ is defined by (\ref{eq:smoothmeasure}). From (\ref{eq:laplace_Z}) and the definition of $\Psi_r(a)$ (just before (\ref{eq:laplace_condition})),
\begin{equation} \label{eq:laplace_jensen}
\mathcal{Z}(r) = \int_{\mathbb{R}^r}\exp\left[-K_{r}(a)-\sum^r_{i=1} \left(w(a_i) + \log\left(m^*_\varepsilon(a_i)\right)\right)\right]\prod^r_{i=1} m^*_\varepsilon(a_i)\hspace{0.02cm}da_i
\end{equation}
Since $\mu^*_\varepsilon(dx) = m^*_\varepsilon(x)\hspace{0.02cm}dx$ is a probability measure, it is possible to apply Jensen's inequality in (\ref{eq:laplace_jensen}). This yields
\begin{align*}
\log\left(\mathcal{Z}(r)\right) 
&\geq - \int_{\mathbb{R}^r} K_r(a)\prod^r_{i=1} m^*_\varepsilon(a_i)\hspace{0.02cm}da_i 
 - \int_{\mathbb{R}^r} \sum^r_{i=1} \left(w(a_i)+\log\left(m^*_\varepsilon(a_i)\right)\right)\prod^r_{i=1}  m^*_\varepsilon(a_i)\hspace{0.02cm}da_i \\[0.1cm]
& = -\int_{\mathbb{R}^r} K_r(a)\prod^r_{i=1} m^*_\varepsilon(a_i)\hspace{0.02cm}da_i 
 - r \int_{\mathbb{R}} \left(w(x)+\log\left(m^*_\varepsilon(x)\right)\right)m^*_\varepsilon(x)\hspace{0.02cm}dx 
\end{align*}
Furthermore, by applying (\ref{eq:multipleK}) to the probability measure $\mu^*_\varepsilon\hspace{0.03cm}$, 
$$
\log\left(\mathcal{Z}(r)\right) \geq - r(r-1)I(\mu^*_\varepsilon)  
 - r \int_{\mathbb{R}} w(x)\hspace{0.03cm}m^*_\varepsilon(x)\hspace{0.02cm}dx 
 - r \int_{\mathbb{R}} \log\left(m^*_\varepsilon(x)\right)\hspace{0.03cm}m^*_\varepsilon(x)\hspace{0.02cm}dx
$$
where both integrals on the right-hand side are finite, as can be seen by noting that $w$ is locally integrable and that $m^*_\varepsilon$ is positive, continuous, and compactly supported. Now, dividing both sides by $r^2$,
\begin{align*}
\frac{1}{r^2}\hspace{0.02cm}\log\left(\mathcal{Z}(r)\right) &\geq \left(\frac{1}{r}-1\right)I(\mu^*_\varepsilon) 
- \frac{1}{r} \int_{\mathbb{R}} w(x)\hspace{0.03cm}m^*_\varepsilon(x)\hspace{0.02cm}dx 
- \frac{1}{r} \int_{\mathbb{R}} \log\left(m^*_\varepsilon(x)\right)\hspace{0.03cm}m^*_\varepsilon(x)\hspace{0.02cm}dx 
\\[0.1cm]
& = \left(\frac{1}{r}-1\right)I(\mu^*_\varepsilon) - \frac{c_\varepsilon}{r}
\end{align*}
for some real number $c_\varepsilon\hspace{0.03cm}$. However, by Lemma \ref{lem:stability}, $I(\mu^*_\varepsilon) \rightarrow I(\mu^*)$ as $\varepsilon \rightarrow 0$, so $I(\mu^*_\varepsilon) \leq I(\mu^*) + \varepsilon^\prime$, for any $\varepsilon^\prime > 0$, whenever $\varepsilon$ is small enough. Replacing this into the last inequality, 
$$
\frac{1}{r^2}\hspace{0.02cm}\log\left(\mathcal{Z}(r)\right) \geq 
\left(\frac{1}{r}-1\right)(I(\mu^*)+\varepsilon^\prime) - \frac{c_\varepsilon}{r}
$$
Then, letting $r \rightarrow \infty$ while $\varepsilon$ and $\varepsilon^\prime$ are fixed, 
$$
 \liminf\, \frac{1}{r^2}\hspace{0.02cm}\log\!\left(\mathcal{Z}(r)\right) \geq -(I(\mu^*)+\varepsilon^\prime)
$$
Finally, since $I(\mu^*) = E$ and $\varepsilon^\prime$ is arbitrary, (\ref{eq:proof_laplace_2}) becomes immediate. 

\subsection{Proofs for Section \ref{sec:cones}} \label{app:cones}

\noindent \textbf{Proof of (\ref{eq:Zcone_asymp})\,:} Dropping the factor $\omega_\beta/r!$ from (\ref{eq:Zcone}),
\begin{equation} \label{eq:Zacone}
Z_\mathfrak{a}(\sigma) = \int_{\mathbb{R}^r} e^{-\Vert a\Vert^2/2\sigma^2}\hspace{0.02cm}\prod_{i<j} |\sinh(a_i-a_j)|^\beta\hspace{0.02cm}da
\end{equation}
The first step is to bring this into the form (\ref{eq:IA}). Let $\mathfrak{a}$ be the hyperplane in $\mathbb{R}^r$ defined by the equation $\sum a_i = 0$. Then,  write $a_i = r^{-1/2}\hspace{0.03cm}x + \bar{a}_i$ where $x \in \mathbb{R}$ and $\bar{a} \in \mathfrak{a}$. Integrating with respect to the new variables $x$ and $\bar{a}$,
$$
 Z_\mathfrak{a}(\sigma) = (2\pi\sigma^2)^{1/2}\hspace{0.03cm}\int_{\mathfrak{a}}
e^{-\Vert \bar{a}\Vert^2/2\sigma^2}\hspace{0.02cm}\prod_{i<j} |\sinh(\bar{a}_i-\bar{a}_j)|^\beta\hspace{0.02cm}d\bar{a}
$$
This will be written, $Z_\mathfrak{a}(\sigma) = (2\pi\sigma^2)^{1/2}\hspace{0.02cm}\bar{Z}_\mathfrak{a}(\sigma)$. The integral $\bar{Z}_\mathfrak{a}(\sigma)$ is of the form (\ref{eq:IA}), with $r$ instead of $r+1$. Putting $\sigma^2 = t/r$ and taking logarithms,
$$
\frac{1}{r^2}\log\!\left(Z_\mathfrak{a}(\sigma)\right) = \frac{1}{2r^2}\log\!\left(2\pi\hspace{0.02cm}t/r\right) + \frac{1}{r^2}\log\!\left(\bar{Z}_\mathfrak{a}(\sigma)\right)
$$
The first term vanishes when $r \rightarrow \infty$, while, by Proposition \ref{prop:equilibrium}, the second term converges to $-E_A = -E_\beta(t)$. Indeed, (\ref{eq:lima}) of Proposition \ref{prop:equilibrium} implies the limit of $(r-1)^{-2}\log\!\left(\bar{Z}_\mathfrak{a}(\sigma)\right)$ is equal to $-E_A\hspace{0.03cm}$. To prove (\ref{eq:Zcone_asymp}), it remains to show that $E_\beta(t) = (\beta/2)\hspace{0.02cm}E_2(\beta t/2)$. To do so, recall that $E_{A} = \inf I_{A}(\mu)$, as stated just before Proposition \ref{prop:equilibrium}. The functional $I_A$ is given by (\ref{eq:functionala}), which can be written
$$
I_{A}(\mu) = \frac{\beta}{2}\left(\frac{c}{\beta t} \int\, x^2\hspace{0.02cm}\mu(dx) - \int\int\,\ell(x-y)\hspace{0.02cm}\mu(dx)\mu(dy)\right)
$$
However, the expression in parentheses is just $I_A(\mu)$ when $\beta = 2$ and $t$ is replaced by $\beta t/2$. Taking the infimum of the left hand side, one obtains $E_A = E_\beta(t)$ (notice that this is the definition of $E_\beta(t)$). For the right-hand side, the infimum is $(\beta/2)$ times the infimum of the expression in parentheses, or $(\beta/2)\hspace{0.02cm}E_2(\beta t/2)$. \\[0.1cm]
\textbf{Proof of (\ref{eq:zcone_asymp})\,:} The starting point is (\ref{eq:zcone_mehta}). Putting $\sigma^2 = t/r$ and taking logarithms, it is easy to show
\begin{equation} \label{eq:z_asymp_proof_cone}
\frac{1}{r^2}\log\!\left(z_\mathfrak{a}(\sigma)\right) = \frac{\beta}{4}\log(t/4)-\frac{\beta}{4}\log(r) + \frac{1}{r^2} \sum^r_{j=1} \log\left(\Gamma(1+j\beta/2)\right) + O\left(\frac{\log(r)}{r}\right)
\end{equation}
Now, note that the sum may be replaced by an integral,
$$
\sum^r_{j=1} \log\left(\Gamma(1+j\beta/2)\right) = \int^r_0 \log\left(\Gamma(1+u\beta/2)\right)du + R(r)
$$
with the error term,
$$
|R(r)| \leq (r\beta/2) \times \sup_{u \in [0,r]} |\psi(1+u\beta/2)|
$$
Here, the $\psi$ function is the logarithmic derivative of the Gamma function~\cite{whittaker}. For large $u$, this has the expansion $\psi(u) = \log(u) - (1/2u) + O(1/u^2)$~\cite{whittaker} (Page 251), which implies that $R(r) = O(r\log(r))$. With this estimate, (\ref{eq:z_asymp_proof_cone}) becomes 
\begin{equation} \label{eq:z_asymp_proof_cone_1}
\frac{1}{r^2}\log\!\left(z_\mathfrak{a}(\sigma)\right) = \frac{\beta}{4}\log(t/4)-\frac{\beta}{4}\log(r) + \frac{1}{r^2} \int^r_0 \log\left(\Gamma(1+u\beta/2)\right)du + O\left(\frac{\log(r)}{r}\right)
\end{equation} 
This will yield (\ref{eq:zcone_asymp}) if one can show
\begin{equation} \label{eq:gammaverage}
  \frac{1}{r^2} \int^r_0 \log\left(\Gamma(1+u\beta/2)\right)du = \frac{\beta}{2}\left(-\frac{3}{4} + \frac{1}{2}\log(r\beta/2)\right) + O\left(\frac{\log(r)}{r}\right)
\end{equation}
Indeed, (\ref{eq:zcone_asymp}) follows immediately by replacing (\ref{eq:gammaverage}) into (\ref{eq:z_asymp_proof_cone_1}) and taking the limit $r \rightarrow \infty$. On the other hand, (\ref{eq:gammaverage}) can be obtained by expressing the integral in terms of the Barnes G-function~\cite{whittaker} (Page 264), and then by using the asymptotic expansion of the Barnes G-function~\cite{barnes}. \\[0.1cm]
\noindent \textbf{Remark\,:} the integrals $Z(\sigma)$ and $z(\sigma)$ in (\ref{eq:Zcone}) and (\ref{eq:zcone}) are not exactly the same as (\ref{eq:Z}) and (\ref{eq:z}), since the domain of integration is $\mathbb{R}^r$ without the restriction $\sum a_i = 0$. As seen in the proof of (\ref{eq:Zcone_asymp}), this restriction merely divides $Z(\sigma)$ by a factor $(2\pi\sigma^2)$.  
In fact, by an identical reasoning, the same is true for $z(\sigma)$. Therefore, the ratio $Z(\sigma)/z(\sigma)$ is~the same, whether the domain of integration is the whole of $\mathbb{R}^r$ or just $\mathfrak{a}$ (\textit{i.e.} $\mathbb{R}^r$ with the above-said restriction). This justifies applying Proposition \ref{prop:free} to this ratio, when $\beta = 2$. The corresponding limit is indeed the one in (\ref{eq:asymptotics_A}) ($A_{r}$ root system). 

\appendices

\section{Perturbative expansion of $Z(\sigma)$} \label{app:perturbation}
The perturbative expansion of $Z(\sigma)$, briefly mentioned in the introduction, was abandoned in the present work. 
This appendix clarifies its underlying idea.

Returning to the proof of Proposition \ref{prop:Ztaylor}, recall (\ref{eq:prooftaylor1}) and (\ref{eq:prooftaylor2}), and note that 
$\mathrm{ad}^2(v)$ is the radial curvature operator,
$R_v : T_{x_o}M \rightarrow T_{x_o}M$,  which is given by~\cite{helgason} (Page 215)
$$
\mathrm{ad}^2(v)\hspace{0.02cm}u = R_v(u) = R(v,u)\hspace{0.02cm}v
$$
with $R$ the Riemann curvature tensor. The determinant (\ref{eq:prooftaylor2}) is equal to 
\begin{align*}
\det\left(\sum^\infty_{n=0} \frac{R^{n}_v}{(2n+1)!} \right)
\end{align*}
As a function of $v$, this determinant has a power series expansion whose non-constant terms are linear combinations of the traces of exterior powers $\mathrm{tr}(\wedge^k\hspace{0.02cm}R^n_v)$ ($k=1,\ldots,d$ and $n = 1, 2,\ldots$). Each term is then a polynomial function of $v$ with coefficients determined by the Riemann curvature tensor. With~this in mind, it is possible to re-write (\ref{eq:prooftaylor1}),
$$
Z(\sigma) = \hspace{0.03cm}\int_{T_{x_o}M}e^{-\Vert v\Vert^2/2\sigma^2}\left(1 + P_{\hspace{0.02cm} 2}(v) + P_{\hspace{0.02cm} 4}(v) + \ldots\right)dv
$$
where each polynomial $P_{\hspace{0.02cm} 2n}$ has degree $2n$ and is moreover $K$-invariant (\textit{i.e.} isotropy-invariant), because $R_v$ is~obviously $K$-invariant. It is now easy to obtain
$$
Z(\sigma) = z(\sigma)\left(1+ f_2\hspace{0.03cm}\sigma^2 + f_4\hspace{0.03cm}\sigma^4 + \ldots\hspace{0.03cm} \right)
$$
where the general coefficient $f_{2n}$ reads
$$
f_{2n} = \frac{1}{(2\pi)^{d/2}}\int_{T_{x_o}M}P_{\hspace{0.02cm}2n}(v)\hspace{0.02cm}e^{-\Vert v\Vert^2/2}\hspace{0.02cm}dv
$$
just as in the introduction. 

For example, the first two coefficients can be found as follows. For $n = 1$, (in Einstein notation)
$$
P_{\hspace{0.02cm}2}(v) = \frac{1}{3!}\mathrm{tr}(R_v) = \frac{1}{3!}g^{kl}\hspace{0.02cm}R_{ikjl}\hspace{0.02cm}v^iv^j
$$
where $g$ is the metric tensor. Applying Wick's theorem merely amounts to adding up the contractions with respect to all
pairs of indices carried by $v$.
Thus, $f_{\hspace{0.02cm} 2} = (1/6)\hspace{0.02cm}g^{ij}g^{kl}\hspace{0.02cm}R_{ikjl}$ and this is $1/6$ of the scalar curvature of $M$. For~$n = 2$, 
$$
P_{\hspace{0.02cm}4}(v) = \frac{1}{5!}\mathrm{tr}(R^{\hspace{0.02cm}2}_v) + \frac{1}{(3!)^2}\mathrm{tr}(\wedge^2\hspace{0.02cm}R^{\phantom 2}_v)
$$
upon expressing this in terms of the Riemann tensor and applying Wick's theorem,
\begin{align*}
f_{\hspace{0.02cm} 4} & = \frac{1}{5!}\left(g^{ii^\prime}\!g^{kk^\prime}+g^{ik}\!g^{i^\prime\!\hspace{0.02cm}k^\prime} + 
g^{ik^\prime}\!g^{ki^\prime}\right)g^{jl^\prime}g^{lj^\prime} R_{ijkl}R_{i^\prime j^\prime k^\prime l^\prime} \\[0.1cm]
                                   & \!\!\!\!\!\!\!+ \frac{1}{2(3!)^2} 
\left(g^{ii^\prime}\!g^{kk^\prime}+g^{ik}\!g^{i^\prime\!\hspace{0.02cm}k^\prime} + 
g^{ik^\prime}\!g^{ki^\prime}\right)\left(g^{j^\prime\!\hspace{0.03cm}l^\prime}g^{jl}-g^{jl^\prime}g^{lj^\prime}\right) R_{ijkl}R_{i^\prime j^\prime k^\prime l^\prime}
\end{align*}
It remains to make sens of such expressions, perhaps through a diagrammatic approach.

\bibliographystyle{IEEEtran}
\bibliography{references}

\end{document}